\documentclass{amsart}

\usepackage[dvipsnames]{xcolor}
\usepackage{tikz}
\usepackage{bm,bbm}
\usepackage{amsmath,amsthm,amssymb}
\usepackage{fullpage}
\usepackage{color}
\usepackage{dsfont}
\usepackage{etoolbox}
\usepackage{comment}
\usepackage{pgf}
\usetikzlibrary{calc}
\usetikzlibrary{patterns}
\usetikzlibrary{arrows}
\usetikzlibrary{decorations.pathreplacing}
\usetikzlibrary{positioning}
\usepackage[utf8]{inputenc}
\usepackage[most]{tcolorbox}
\usepackage{adjustbox}
\usepackage[final]{pdfpages}
\usepackage{blkarray}
\usepackage{wrapfig}
\usepackage{ marvosym }
\usepackage{enumitem}
\usepackage{hyperref}
\theoremstyle{plain}
\newtheorem{theorem}{Theorem}
\newtheorem{proposition}[theorem]{Proposition}
\newtheorem{corollary}[theorem]{Corollary}
\newtheorem{lemma}[theorem]{Lemma}

\theoremstyle{definition}

\theoremstyle{remark}
\newtheorem{remark}[theorem]{Remark}
\newtheorem*{remark*}{Remark}

\newcounter{itemrefcounter}
\makeatletter
\newcommand\myitem[1][]{\item[#1]\refstepcounter{itemrefcounter}\def\@currentlabel{#1}}
\makeatother

\def\E{\mathbb E}

\def\P{\mathbb P}

\def\R{\mathbb R}

\def\mc{\mathcal}
\def\ms{\mathsf}

\def\s{\sigma}
\def\su{\subseteq}

\def\e{\varepsilon}

\def\de{\delta}

\def\one{\mathbbmss{1}}

\def\ff{\infty}

\def\MM{\mc M}

\def\PP{\mc P}
\def\FF{\mc F}
\def\CC{\mc C}

\def\DD{\mc D}

\def\XX{\mc X}

\def\AA{\mc A}

\def\f{\frac}

\def\im{\item}

\def\ti{\times}

\def\supp{\ms{supp}}

\renewcommand\le{\leqslant}
\renewcommand\ge{\geqslant}

\def\X{\mathbb X}

\def\co{\colon}

\def\N{\mathbb N}

\usepackage[foot]{amsaddr}

\begin{document}

\title{Large deviations of one-hidden-layer neural networks}

\date{\today}

\author{Christian Hirsch$^1$}
\address{$^1$Aarhus university}
\email{hirsch@math.au.dk}

\author{Daniel Willhalm$^2$}
\address{$^2$University of Groningen}
\email{d.willhalm@rug.nl}

\begin{abstract}
    We study large deviations in the context of stochastic gradient descent for one-hidden-layer neural networks with quadratic loss. We derive a quenched large deviation principle, where we condition on an initial weight measure, and an annealed large deviation principle for the empirical weight evolution during training when letting the number of neurons and the number of training iterations simultaneously tend to infinity. The weight evolution is treated as an interacting dynamic particle system. The distinctive aspect compared to prior work on interacting particle systems lies in the discrete particle updates, simultaneously with a growing number of particles.
\end{abstract}

\maketitle

\section{Introduction}
\label{int_sec}

Large artificial neural networks (ANN) are nowadays indispensable in industry. A milestone was the convolutional neural network AlexNet in 2012, a model for image classification that won the ImageNet Large Scale Visual Recognition Challenge and outperformed all other models by a large margin. Also the use of graphics processing units to train the model via stochastic gradient descent (SGD) was a breakthrough in the field. This allowed the size of neural networks to grow significantly since then, and they became state of the art in many fields. The downside is that AlexNet with its over 9000 neurons also marks the point at which our understanding of the behavior of individual neurons and their collective interactions fades, as the network's operations are often obscured by the sheer volume of computations occurring simultaneously. But shedding light on the large-scale neural networks is of great interest. This concerns not only for deep neural networks such as AlexNet which had eight hidden layers, three of which were fully connected ones, but also shallow neural networks with only one hidden layer. From a theoretical point of view, the exploration of whether the intricate training process of neural networks conforms to fundamental principles of probability theory, such as laws of large numbers (LLN), central limit theorems or large deviation principles (LDP) has become a vibrant research field.
More precisely, in this paper, we focus on supervised learning using a neural network consisting of one input layer, one hidden layer and an output layer trained via SGD with the goal of minimizing the square loss against empirical data. We study the limit of a growing number of nodes in the hidden layer and a growing number of training iterations. The neural network is structured with weights and an activation function, and the aim of training is to adjust these weights to reduce the squared error between the predicted output and the true label of a data point drawn from a given distribution. The weights are drawn at random from an initial weight distribution. When the activation function is nonlinear, this finite-dimensional optimization problem becomes nonconvex, which means finding the optimal solution using any form of gradient descent is not straightforward.
However, by reformulating the problem in a larger mathematical framework, particularly by lifting it to a space of probability distributions over the network parameters, the originally nonconvex problem can be turned into a convex one. This transformation that the reader is guided through in Section \ref{sec:mod}, opens up the possibility of finding a global minimum more reliably using gradient methods within this larger space.
This leads to several key questions, which are central to the paper and related literature. We summarize them as follows:
\begin{itemize}[leftmargin=0.7cm]
\item[Q1] How does SGD applied to the original finite-dimensional problem relate to the gradient descent applied in the larger space of measures when the number of neurons and the number of training steps grows infinitely large? Essentially, this question explores the transition from the finite-dimensional to the measure-based framework in neural network training.
\item[Q2] If SGD in the measure space leads to convergence of the measures when training time tends to infinity, does the objective function, i.e., the loss also reach its minimum? Further, is there a particular distribution of parameters that SGD converges to as optimal solution?
\item[Q3] When does SGD in the measure space reliably converge to a solution?
\end{itemize}
Regarding Q1, a breakthrough was obtained in \cite{siri}, where an LLN was proven, and \cite{siriCLT}, which illustrates the derivation of a central limit theorem. In particular, \cite{siri} shows that the typical evolution of the weights during SGD training satisfies a McKean-Vlasov equation.
Q2 is confirmed by \cite{chizat2018global}. A general convergence of gradient descent in the measure space as suggested in Q3 has not been possible to show. However, certain modifications such as adding entropic regularization can ensure reliable \cite{meimontanaringuyen2018,kaitong2021} and even exponential \cite{nitandawusuzuki2022,chizat2022} convergence. Other results for SGD for one-hidden-layer ANNs are for instance \cite{rotskoff2022,jungleeleeyang2023,apollonio2023normal,caron2023overparameterised,hirschneumannschmidt}. Finally, besides the LLN,  also the propagation of chaos is a highly relevant topic in this area since it guarantees asymptotic decorrelation of the individual particles for large systems. Such decorrelation properties are particularly interesting if they are quantitative and uniform in time, and we refer the reader to \cite{lacker} for recent results in this direction.

In this work, we add to the understanding of Q1 and deal with the far more delicate problem of characterizing the weight evolution in atypical realizations and establish an LDP for the evolution of the weights during the training of a one-hidden-layer perceptron when number of nodes in the hidden layer and the number of training iterations tend to infinity simultaneously. However, large deviations theory has great potential for practical applications. This is because this theory allows converting a probabilistic question on rare events into a deterministic entropy-optimization problem. That being said, due to the complexity of the involved spaces such a minimization cannot be carried out in closed form so that approximations would need to be used. However, developing such approximations is a difficult task on its own and is therefore outside the scope of the present paper. Moreover, importance sampling is also a useful application of large deviation theory.

The large-deviation analysis of classical SGD systems is one instance of large deviations theory for interacting particle systems, where the weights act as particles evolving over time. While the theory of large deviations for scalar random variables is well understood, dealing with interacting particle systems is far more challenging. Only recently, a general framework was developed in \cite{bufi} which provides a blueprint for addressing such large deviations problems. This approach is based on a weak-convergence approach presented in \cite{dupuisellis}. The idea is to first derive a variational identity for the relevant rare-event probability, and then analyze this expression in the limit. While \cite{bufi} is concerned with interacting particle systems driven by Brownian motions, their arguments admit much more flexibility. This observation is also taken up in the monograph \cite{dupuis}, and we leverage this theory in our investigations.

The difficulty in the analysis of the weight evolution during SGD stems from dealing with a setting that involves discrete updates of the weights while dealing with an increasing number of particles. The existing literature so far only considered each of these scenarios individually but not jointly. For example, \cite{dupuis} investigates large deviations of an increasing number of interacting particles where dynamics are determined by a Brownian motion in continuous time. The theory of \cite{bufi} allows for discrete updates but only for one single particle. Our research combines these two settings. Another difficulty stems from \emph{simultaneous transitions} of the weights. That is, one new datapoint causes an update of \emph{all} weights. Such transitions are difficult to handle in the general framework of \cite{dupuis}, where each particle has its own Brownian motion that is independent of other Brownian motions. For instance, \cite{drw} deals with a discrete setting where a new datapoint can lead to the simultaneous modification of a finite fixed number of weights. In a continuous setting, \cite{budconroy} derive large deviations, where particles are simultaneously affected by a common Brownian motion. However, in our work, we deal with updating all weights simultaneously in a discrete setting. We prove a quenched LDP for the empirical weight evolution measure, conditioning on the initial measure and derive an annealed LDP from it. We give conditions under which the LDPs can be applied. These include boundedness, continuity and differentiability of the activation function, a compactly supported initial weight distribution, the existence of exponential moments related to the data distribution and a uniqueness property of a stochastic differential equation (SDE). We further show that a compactly supported data distribution is a sufficient condition for the uniqueness property and we derive a weak LLN from the LDP.

The weak convergence approach is invoked to prove the quenched LDP, where the work lies in establishing tightness and identifying the limiting distribution. In contrast to \cite{bufi}, who deal with the space of continuous functions with respect to the uniform topology, the objects we investigate are measures on a Skorokhod space. This brings additional difficulties when establishing tightness, where we have to directly deal with compactness in a Skorokhod space and make use of the precise dynamics of the training via SGD. We like to point out that it is also possible to avoid dealing with those technicalities and work with the uniform topology instead because in the limit weight trajectories will be continuous with respect to the uniform topology. One could use this and work with a linear interpolation of the weights instead of a constant interpolation which is what we used. The strategy would then be to show that constant and linear interpolation in this case admit the same large deviations. However, while some arguments in the section on tightness and the section concerned with the upper bound would simplify, the theory of working with the Skorokhod topology is well established in \cite{billingsley} and also \cite{kallenberg} contains helpful theorems when dealing with continuity in the limit. Further, more notation would be necessary to introduce linear interpolation and additional arguments would be required to show that the large deviations indeed transfer. Therefore, we decided to work with the Skorokhod topology.
Also compared to \cite{siri} our proof of tightness is more involved, as they only have to consider a fixed data distribution and not random ones that arise due to the representation. Our identification of the limit is associated with \cite{bufi} but we do not need to deal with a local martingale problem. However, our Skorokhod setting requires a more careful analysis, and a lot of our work involves dealing with truncations, which \cite{bufi} did not have to consider. The annealed LDP is derived from the quenched LDP by using a general result that can be found in \cite{annealedLDP} and \cite{Ganesh}. The proof that a compactly supported data distribution is a sufficient condition for the uniqueness property is based on \cite{siri}, who proved a very similar result but with the mentioned difference that stems from the representation, we have to deal with more measures than only the fixed data distribution.

A possible extension to multi-hidden-layer perceptrons is still very far-fetched. Quite recently, \cite{siriDeep} started to tackle limit theory for deep neural networks and in \cite[Section 4.3]{siriDeep} they explain the challenges that arose when trying to extend the LLN or central limit theorem from the single-hidden-layer approach to a multi-hidden-layer neural network. Further attempts have been made in \cite{nguyenphan2023} and \cite{araujo2019meanfield}.

The article is structured as follows. In Section \ref{sec:mod}, we introduce the model and the main results. In Section \ref{sec:proof_overview_quenched}, we provide the proof of the quenched LDP and state all lemmas necessary for it. Section \ref{sec:proof_overview_annealed} deals with the lemmas for and the proof of the annealed LDP.
Section \ref{sec:proof_overview_quenched} is devoted to the detailed proofs of the lemmas that were announced in Section \ref{sec:proof_overview}.
In Section \ref{sec:tight}, we prove tightness of the respective measures. In Section \ref{sec:ident}, we identify the limiting distribution of the empirical weight evolution measure and prove goodness of the rate function.
Finally, in Section \ref{sec:unq}, we prove a sufficient condition for the uniqueness property in the annealed and quenched LDP theorems and derive a weak LLN from the previous results.

\section{Model and main results}
\label{sec:mod}
To present the model, we follow in broad strokes the notation of \cite{siri}.

\subsection{Network architecture}
We investigate a network consisting of an input layer with $d' \ge 1$ nodes and one hidden layer of $n \ge 1$ nodes. The output layer is linear and is connected to the hidden layer by weights $c^{1, n}, \dots, c^{n, n} \in \R$. The hidden layer is connected to the input layer through an activation function $\s$ and weight vectors $w^{1, n}, \dots, w^{n, n}\in \R^{d'}$. Taken together, the weights $( \theta_0^{i,n})_{i \le n} := (c^{i,n}, w^{i,n})_{i \le n}\su \R^d := \R^{d' + 1}$ form the initial network parameters.  We assume that the initial weights are i.i.d.\  random variables drawn from a distribution $\nu$ on $\R^d$. Then, we can represent the output of the ANN with input $z \in \R^{d'}$  by  
\begin{align}
	\label{eq_g}
	F(z,  \theta_0^n) = \int_{\R^d} c \s(z^\top w)  \theta_0^n({\rm d}  \theta),
\end{align}
where  $ \theta_0^n := \f1n\sum_{i \le n} \delta_{ \theta_0^{i, n}}$ is the empirical initial weight measure. We also assume that the training data is of the form $ X:= (Z, Y) $, where $Z$ represents the input vector and $Y$ the output value, following some distribution $\pi$ on $\R^d$. Throughout the paper we denote by $\|\cdot\|$ the Euclidean norm on $\R^d$ or $\R^{d'}$ and by $\|\cdot\|_\ff$ the supremum norm for real valued functions.

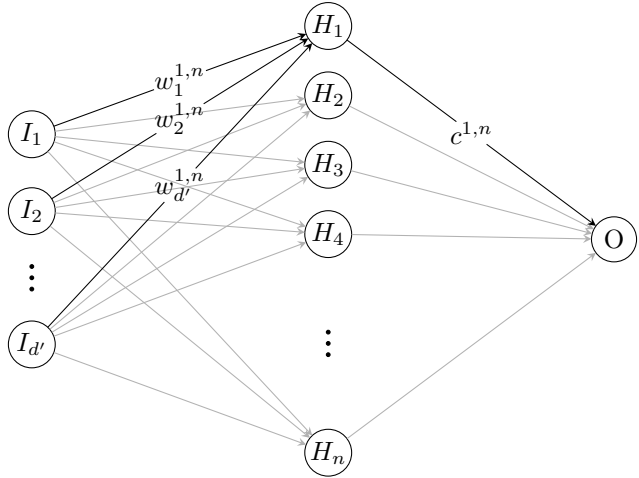
\begin{figure}[ht]
\centering
\begin{tikzpicture}[
    neuron/.style={circle, draw, minimum size=0.55cm, text width=0.55cm, inner sep=0, outer sep=0, align=center},
    dot/.style={circle, draw, fill, inner sep=0, outer sep=0, minimum size=0.5mm},
    weight/.style={sloped, anchor=south, font=\footnotesize, auto=false, inner sep=1pt}
]

\node[neuron] (input1) {$I_1$};
\node[neuron, below=0.4cm of input1] (input2) {$I_2$};
\node[dot, below=0.4cm of input2] (inputdot1) {};
\node[dot, below=0.1cm of inputdot1] (inputdot2) {};
\node[dot, below=0.1cm of inputdot2] (inputdot3) {};
\node[neuron, below=0.4cm of inputdot3] (input3) {$I_{d'}$};

\node[neuron, above right=1cm and 3.5cm of input1] (hidden1) {$H_1$};
\node[neuron, below=0.3cm of hidden1] (hidden2) {$H_2$};
\node[neuron, below=0.3cm of hidden2] (hidden3) {$H_3$};
\node[neuron, below=0.3cm of hidden3] (hidden4) {$H_4$};
\node[neuron, below right=1cm and 3.5cm of input3] (hidden5) {$H_n$};
\coordinate (middledots) at ($(hidden4)!0.5!(hidden5)$);
\node[dot] (hiddendot1) at ($(hidden4)!0.5!(hidden5)$) {};
\node[dot, above=0.1cm of hiddendot1] (hiddendot2) {};
\node[dot, below=0.1cm of hiddendot1] (hiddendot3) {};

\coordinate (hiddenmiddle) at ($(hidden1)!0.5!(hidden5)$);
\node[neuron, right=3.5cm of hiddenmiddle] (output) {O};

\foreach \i in {2,3,4,5}
    \draw[-stealth, black!30] (hidden\i) -- (output);
    
\foreach \i in {1,2,3}
    \foreach \j in {2,3,4,5}
        \draw[-stealth, black!30] (input\i) -- (hidden\j);

\path (input1) -- (hidden1) node[pos=0.5, inner sep=0] (node1){$w_1^{1,n}$}[behind path];

\path (input2) -- (hidden1) node[pos=0.5, inner sep=-0.6pt] (node2){$w_2^{1,n}$}[behind path];
        
\path (input3) -- (hidden1) node[pos=0.5, inner sep=0] (node3){$w_{d'}^{1,n}$}[behind path];

\path (hidden1) -- (output)
          node[pos=0.5, inner sep=0] (node4){$c^{1,n}$}[behind path];

\foreach \i in {1,2}
    \foreach \j in {1}
    	\path (input\i) edge[-stealth, to path=--(node\i)--(\tikztotarget)] (hidden\j);
    	
\path (input3) edge[-stealth, to path=--(node3)--(\tikztotarget)] (hidden1);

\path (hidden1) edge[-stealth, to path=--(node4)--(\tikztotarget)] (output);
\end{tikzpicture}
\caption{Architecture of the considered one-hidden-layer neural network}
\end{figure}

\subsection{Stochastic gradient descent}\label{sec:sgd}
The next step is to implement SGD for the quadratic loss function during a training interval $[0, T]$, where $T>0$ is fixed. We assume that weights are updated every $1/n$ units of time.
In the training step associated with $k\le n T$ the weight vectors for every node are changed using only one incoming data point. More precisely, we let 
$$(X_{k, n})_{k \le nT} = \big((Z_{k, n},Y_{k, n} )\big)_{k \le nT}$$
be an i.i.d.\ family of $\pi$-distributed random vectors serving as input for the weight updates.
The objective is to minimize the quadratic loss of the network, which is given by
\begin{equation}\label{eq:objective_finite_dim}
	\min_{(c^{j,n},w^{j,n})_{j\le n}\subseteq\R^d} \int_{\R^d} \frac12 \Big(y - \frac1n \sum_{i\le n} c^{i,n} \s(z^\top w^{i,n})\Big)^2 \pi({\rm d} (z,y)),
\end{equation}
similar to \cite[(1.2)]{siri}. When training the network via SGD we are interested in how small weight changes affect the loss function for a single training data point $(z,y)\in\supp(\pi)$. This is quantified by the gradient of the loss function 
\begin{equation}\label{eq:gradient_loss}
	\nabla_{(c^{j,n},w^{j,n})_{j\le n}}\frac12 \Big(y - \frac1n \sum_{i\le n} c^{i,n} \s(z^\top w^{i,n})\Big)^2.
\end{equation}
The $j$th component of the gradient (in the sense of the weights attached to the $j$th neuron in the hidden layer) is then given by
\begin{equation}\label{eq:gradient_loss_calculated}
	\bigg(\Big(y - \frac1n \sum_{i\le n} c^{i,n} \s(z^\top w^{i,n})\Big)\s(z^\top w^{j,n}), \Big(y - \frac1n \sum_{i\le n} c^{i,n} \s(z^\top w^{i,n})\Big)c^{j,n}\s'(z^\top w^{j,n})z\bigg),
\end{equation}
which is to be understood as a $d$-dimensional vector. As mentioned in the introduction, this optimization problem in \eqref{eq:objective_finite_dim} is nonconvex if the activation function $\s$ is nonlinear. However, lifting  to the space of probability measures on $\R^d$ gives rise to a convex optimization problem
$$\min_{\theta\in\PP(\R^d)} \int_{\R^d} \frac12 \big(y - F(z, \theta)\big)^2 \pi({\rm d} (z,y)),$$
where $\PP(\R^d)$ denotes the set of all probability measures on $\R^d$. The gradient flow for the objective in the broader space of measures is given by a Wasserstein gradient flow, which is explained in detail in \cite{kaitong2021} or in \cite[Remark 1.5]{siri}. We remind the reader that our goal is to investigate how SGD with respect to the finite-dimensional objective in the limit of neurons in the hidden layer and number of training steps relates to the Wasserstein gradient flow in the infinite-dimensional space of probability measures in terms of large deviations.
To that order we would like to keep track of the whole training process and we henceforth investigate the empirical measure that tracks the trajectories of the evolving weights
\begin{equation}\label{eq_emp_measure}
 \theta^n := \f1n\sum_{i\le n} \de_{( \theta_t^{i, n})_t},
\end{equation}
where $( \theta_t^{i, n})_t$ is interpreted as a function that maps a time point $t$ to $\R^d$, which is constructed as follows.

Using the training data $X_{k, n}$ at training step $k\le nT$, we perform an SGD update to all parameters $( \theta_{(k-1)/n}^{i, n})_{i\le n}$ and store the new weights in $(\theta_{k/n}^{i,n})_{i\le n}$. In the first step, we start with the unaltered parameter family $( \theta_0^{i, n})_{i\le n}$. Formally, this means for every $i\le n$ the weight $ \theta_{k/n}^{i,n}$ is characterized via the update equation
\begin{equation}
	\label{eq:c2}
	 \theta^{i, n}_{k/n} :=  \theta^{i, n}_{(k-1)/n} + \e_n A(X_{k, n},  \theta^{i, n}_{(k-1)/n};  \theta_{(k-1)/n}^n), 
\end{equation}
where $\e_n:=1/n$ is the learning rate, $ \theta_{(k-1)/n}^n := \f1n\sum_{i \le n} \delta_{ \theta_{(k-1)/n}^{i, n}}$ and where the gradient is given by \eqref{eq:gradient_loss} and \eqref{eq:gradient_loss_calculated} or more precisely for $x= (z, y)$,
\begin{align}\label{eq:grad}
	\begin{split}
	A\big(x, (c, w) ;  \theta_{(k-1)/n}^n\big) 
	&:= \nabla_{(c,w)} \frac12\big(y - F(z, \theta_{(k-1)/n}^n)\big)^2 \\
	&\;= \big(g(x,  \theta_{(k-1)/n}^n) \s(z^\top w),g(x,  \theta_{(k-1)/n}^n)c\s'(z^\top w)z \big),
	\end{split}
\end{align}
and  $g\big((z, y),  \theta^n_{(k-1)/n}\big) := y - F(z,  \theta^n_{(k-1)/n})$.
There are in total $\lfloor nT\rfloor$ SGD steps, updating each of the $n$ particles using $\lfloor nT\rfloor$ iid training samples. Assuming that the parameters remain unchanged between times $(k-1)/n$ and $k/n$, we can define
$$ \theta_t^{i,n} := \theta_{\lfloor nt\rfloor/n}^{i,n}$$
for every $t\in[0, T]$ and $i\le n$, which defines $ \theta^n$ as a random element in the set $\PP(\XX)$ of all probability measures on the Skorokhod space $\XX:= D([0, T], \R^d)$ of càdlàg trajectories. The goal of this work is to derive an LDP for $ \theta^n$. We note that from the space of probability measures on trajectories in $\XX$, it is possible to specialize to the setting of evolutions of measures on $\R^d$ by projecting to the corresponding time marginals. Therefore our setting is (slightly) broader than the one from \cite{siri}.

Before stating the large deviations asymptotics, we discuss the LLN for $ \theta^n$. For each $t\in[0, T]$, the empirical weight measure $ \theta_t^n$ is identical to the marginal distribution of $ \theta^n$ at time $t$. Then, loosely speaking, for large $n$, the discrete updates through $\pi$-distributed random vectors are replaced by an integration with respect to ${\rm d} t\otimes\pi({\rm d} x)$. Then, the limiting distribution $\eta_t$ at any time $0\le t\le T$ is characterized by the solution to the McKean-Vlasov equation
\begin{align}
	\label{eq:ev}
	\theta_t = \theta_0 +\int_0^t \int_{\R^d} A(x, \theta_s; \ms{Law}(\theta_s)) \pi({\rm d} x) {\rm d} s,\qquad \theta_0\sim \nu,
\end{align}
which is essentially the main result of \cite{siri}. Note that the integral in \eqref{eq:ev} is to be understood as $d$-dimensional vector of integrals. Further, we point out that \eqref{eq:ev} is a stochastic representation of a Wasserstein gradient and refer to \cite{kaitong2021} and \cite[Remark 1.5]{siri} again for more details.

While the limiting evolution equation \eqref{eq:ev} describes the SGD behavior for typical realizations of $ \theta^n$, the contribution of our work will leverage large deviations theory for describing atypical realizations of $ \theta^n$. Loosely speaking, the main idea is that large deviations of the sequence of empirical weight evolution measures $( \theta^n)_n$ are caused by atypical realizations of training data $(X_{k, n})_{k,n}$ and of the weight initialization. We can imagine that those data points and the initial weight configuration are drawn from tilted distributions of the training data and the initial weights, i.e., by measures that are absolutely continuous with respect to $\pi$ and $\nu$, respectively. The tilting of the data points can be described by a measure $\rho$ in 
$$\MM := \MM_\ms{fin}(\X_T)$$
of finite Borel measures on $\X_T := [0, T]\ti \R^d$ such that $\rho$ is the product measure of the Lebesgue measure on $[0,T]$ and a stochastic kernel on $\R^d$ given $[0,T]$ that we denote by $\rho_t := \rho({\rm d} x\mid t)$. We also define $\pi_T$, a product measure on $\X_T$ given by ${\rm d} t\otimes\pi({\rm d} x)$. Next, writing $H(\,\cdot \mid \cdot\,)$ for the relative entropy, we set
$$R(\rho) := H(\rho \mid \pi_T) = \begin{cases} \int_{\X_T} \log\big(\f{{\rm d} \rho}{{\rm d} \pi_T}\big) \rho({\rm d} (t,x)), & \text{if } \rho \ll \pi_T,\\ \ff, & \text{otherwise}\end{cases}.$$
Further, for a probability measure $\eta$ on $\XX$, $\rho\in\MM$ and a stochastic process $( \theta_t)_t$ with distribution $\eta$, we define the evolution equation
\begin{equation}\label{eq:evt}
	 \theta_t -  \theta_0 = \int_0^t \int_{\R^d} A(x, \theta_s; \ms{Law}(\theta_s)) \rho_s ({\rm d} x) {\rm d} s, \qquad  \theta_0 \sim \nu.
\end{equation}

We then assume the following properties.
\begin{enumerate}[leftmargin=\widthof{\textbf{(WCOMP)}}+\labelsep]
\myitem[\textbf{(CONT)}] \label{CON} The activation function $\s$ is Lipschitz continuous and differentiable with $\s'$ also being Lipschitz continuous. Denote a Lipschitz constant for both cases by $L_\s > 1$. Additionally, $\s$ is bounded by some $C_\s\ge 1$.
		\myitem[\textbf{(DEXP)}] \label{DCOMP} Suitable moment generating functions related $\pi$ exist at every point, i.e., for all $c>0$ it holds that $\E[e^{c \|(1, Z)\|^{2}}] \vee \E[e^{c |Y|^{4}}] < \ff$.
\myitem[\textbf{(UNQ)}] \label{UNQ} For every $ \rho\in\MM$ with $R( \rho)<\ff$, the evolution equation \eqref{eq:evt} corresponds to a unique weak solution.  The remark below gives some explanation about what is meant by a correspondence to a weak solution.
\myitem[\textbf{(WCOMP)}] \label{WCOMP} The initial weight distribution $\nu$ has compact support, i.e., there exists $C_\nu \ge 1$ such that $\supp(\nu) \subseteq B_{C_\nu} (0)$, the closed ball in $\R^d$ with radius $C_\nu$ centered at $0$.
\end{enumerate}

\begin{remark}
	That \eqref{eq:evt} corresponds to a unique weak solution means that there is a unique distribution $\eta$ such that $( \theta_t)_t$ with distribution $\eta$ and filtration generated by $( \theta_t)_t$ is a weak solution to \eqref{eq:evt}. 

We stress that the model considered in \cite{bufi} involves an interacting particle system driven by Brownian motion, which is why there weak solutions need to be considered in the fully general setting of SDEs presented in \cite{karatzas1991}. In contrast, in this respect our setting is much less technical since the randomness in the evolution equations  comes only from the initial condition.
\end{remark}

Throughout this paper, we use $\PP(\,\cdot\,)$ to denote the set of all probability measures defined on the Borel sets of a given space. That means $ \theta^n$ is a random element in $\PP(\XX)$ for each $n\in\N$.
The key idea for deriving the large-deviation asymptotics is that the tilting of the distribution of the training data and tilting of the initial weight distribution leads to a tilted evolution equation. For this purpose, we let $(\rho, \eta)\in\MM\times\PP(\XX)$, where the measure $\rho$ represents the data point trajectory measure and $\eta$ the weight evolution trajectory distribution. Further, for $t\in[0,T]$, we define the marginal distribution of $\eta$ projected to time $t$ by $\eta_t$. Note that in the Skorokhod topology, such projections are measurable, see \cite[Theorem 12.5]{billingsley}.
The probabilistic costs associated with the tilting are described by the relative entropy $R(\rho)$.

Similarly as in \cite{bufi}, to describe the large-deviation asymptotics, we need to define a family of admissible measures $\PP_\ff^\nu$ on the state space $\X_T \ti \XX$. More precisely, we write $\PP_\ff^\nu$ for the family of all \emph{admissible} measures in $\MM \ti \PP(\XX)$. Then, we say that $(\rho,\eta)\in\PP_\ff^\nu$ if
\begin{enumerate}
\im[1.] the tilted data distribution $\rho$ satisfies $R(\rho) < \ff$. In particular, $\supp(\rho) \subseteq [0,T]\times\supp(\pi)$;
\im[2.] $\eta$ defines a weak solution for the evolution equation in the sense that for a stochastic process $( \theta_t)_t$ with distribution $\eta$ it holds that $\eta$-almost surely, the evolution equation from \eqref{eq:evt} is satisfied, i.e.,
\begin{equation}\label{eq:SDE}
 \theta_t -  \theta_0 = \int_0^t \int_{\R^d} A(x, \theta_s; \ms{Law}(\theta_s)) \rho_s ({\rm d} x) {\rm d} t, \qquad  \theta_0 \sim \nu;
\end{equation}
\im[3.] the distribution $\eta$ shall have support in the set $\CC$ of continuous functions from $[0,T]$ to $\R^d$.
\end{enumerate}

Now, we can state the main result of this article.
\begin{theorem}[Annealed LDP for $ \theta^n$]
\label{thm:annealed}
Assume that \ref{CON}, \ref{DCOMP}, \ref{UNQ} and \ref{WCOMP} are satisfied. Then, the family of empirical measures $( \theta^n)_{n \ge 1}$ satisfies the LDP in $\PP(\XX)$ with respect to the weak topology with rate function 
$$J(\eta) := \inf_{\nu_0\in\PP(\R^d)} \big(H(\nu_0\mid\nu) + I_{\nu_0}(\eta)\big)$$
for $\eta\in\PP(\XX)$, where
$$I_{\nu_0}(\eta) := \inf_{(\rho,\tilde\eta) \in \PP_\ff^{\nu_0}\co  \tilde\eta= \eta} \frac1{T} R(\rho).$$
\end{theorem}

Here, $I_{\nu_0}$ represents the cost of a tilted data distribution given $\nu_0$ as initial weight distribution. In $J$ the costs of a tilted initial weight distribution are added to the rate function. Essentially, an annealed LDP is an LDP for a sequence of marginal distributions given by a product measure that can be derived if the corresponding sequences of the other marginal probability distribution and the conditional
probability distributions satisfy LDPs. The LDP for the conditional probability distribution is called quenched LDP, see \cite{annealedLDP,Ganesh}.

Next, we would like to emphasize a sufficient condition for the uniqueness requirement \ref{UNQ}. In particular, we show that the uniqueness requirement is satisfied if the data distribution $\pi$ has compact support.
\begin{theorem}[Compact data support implies uniqueness]\label{thm:unq}
Assume that \ref{CON} and \ref{WCOMP} are satisfied. If the data distribution $\pi$ has compact support, then \ref{UNQ} is satisfied.
\end{theorem}
Further, we can recover a weak LLN from Theorem \ref{thm:annealed} that gives some additional information on the speed of convergence for the main result from \cite{siri}. For this, $\eta(f)$ denotes the integral of $f\colon\XX\to\R$ with respect to $\eta\in\PP(\XX)$.
\begin{corollary}[Weak LLN for $ \theta^n$]\label{cor:weakLLN}
Assume that \ref{CON}, \ref{DCOMP}, \ref{UNQ} and \ref{WCOMP} are satisfied. Then, for any bounded and continuous real function $f$ on $\XX$ and $\e>0$ it holds that
$$\P( | \theta^n (f) - \eta (f) | \ge\e) \le \exp\Big(-n \inf_{\eta''\in\{\eta'\in\PP(\XX)\co | \eta' (f) - \eta (f) | \ge\e\}} J(\eta'')\Big) \overset{n\to\ff}{\longrightarrow} 0,$$
where $\eta$ is the (weakly unique) solution of \eqref{eq:evt} for $ \rho_s = \pi$.
\end{corollary}
The proofs of Theorem \ref{thm:unq} and Corollary \ref{cor:weakLLN} are given in Section \ref{sec:unq}.

The proof of the annealed LDP relies on a quenched one. To be able to investigate a quenched LDP, instead of investigating $ \theta^n$, we generalize it to some extent and investigate a wider range of distributions on $\PP(\XX)$ that are not necessarily counting measures. As a first step, we allow the initial weight distribution to depend on $n$ and generalize the update equation \eqref{eq:c2} with that. For this, let $\nu_0^n\in\PP(\R^d)$ such that $\nu_0^n$ converges weakly to $\nu$. Assume the following condition similar to \ref{WCOMP}
\begin{enumerate}[leftmargin=\widthof{\textbf{(WCOMP')}}+\labelsep]
	\myitem[\textbf{(WCOMP')}] \label{WCOMP'} The initial weight distributions $\nu_0^n$ have a uniformly compact support, i.e., there exists $C_\nu>0$ such that $\supp(\nu_0^n)\subseteq B_{C_\nu}(0)$ for all $n\in\N$.
\end{enumerate}

Next, we need to construct the weight evolution in a more technically involved way.
\begin{itemize}
	\item[Step 1] For $(c,w), x\in\R^d$ and $\mu\in\PP(\R^d)$ we define the update step
	$$\ms{SGD}^{x,\mu} (c,w) := \ms{SGD}(x, (c,w); \mu) := (c,w) + \e_n A(x, (c,w); \mu).$$
	\item[Step 2] With this define the pushforward of $\nu_0^n$ given $\ms{SGD}^{x,\nu_0^n}$ by $\nu_1^n(x) := \ms{SGD}_*^{x,\nu_0^n} (\nu_0^n)$
	and for $\lfloor nT\rfloor$ data points $x_{1,n},\dots,x_{\lfloor nT\rfloor,n}\in\R^d$ we proceed by recursively setting
	$$\nu_k^n((x_{i,n})_{i\le k}) := \ms{SGD}_*^{ x_{k,n}, \nu_{k-1}^n((x_{i,n})_{i\le k-1})} (\nu_{k-1}^n((x_{i,n})_{i\le k-1}))$$
	for $k\le nT$. Recall that for a Borel set $B\subseteq\R^d$, the probability measure $\ms{SGD}_*^{x,\mu}(\mu)(B)$ is given by $\ms{SGD}_*^{x,\mu}(\mu)(B) := \mu((\ms{SGD}^{x,\mu})^{-1}(B))$.
	\item[Step 3] For $\tilde \theta_0^n\sim\nu_0^n$, we let $\tilde \theta_{1/n}^n(x_{1,n}) := \tilde \theta_0^n + \e_n A(x_{1,n}, \tilde \theta_0^n; \nu_0^n)$ and for $k\le nT$ we recursively set
	\begin{equation}\label{eq:evt2}
		\tilde \theta_{k/n}^n((x_{i,n})_{i\le k}) :=\ms{SGD}(x_{k,n}, \tilde \theta_{(k-1)/n}^n((x_{i,n})_{i\le k-1}), \nu_{k-1}^n((x_{i,n})_{i\le k-1})).
	\end{equation}
	This gives an element $(\tilde \theta_t^n((x_{i,n})_{i\le \lfloor nt\rfloor}))_t$ in $\XX$ by setting
	$$\tilde \theta_t^n((x_{i,n})_{i\le \lfloor nt\rfloor}) := \tilde \theta_{\lfloor nt\rfloor/n}^n((x_{i,n})_{i\le \lfloor nt\rfloor})$$
	for $t\in[0,T]$.
	\item[Step 4] We define a deterministic element in $\PP(\XX)$ by
	$$\eta^n((x_{i,n})_{i\le \lfloor nT\rfloor}) := \ms{Law}_{\nu_0^n}((\tilde \theta_t^n((x_{i,n})_{i\le \lfloor nt\rfloor}))_t),$$
	where $\ms{Law}_{\mu}$ denotes the law of the argument with respect to the initial distribution $\mu$.\\
	If we instead of fixed data points $x_{1,n},\dots,x_{\lfloor nT\rfloor,n}$ consider the random data points $X_{1,n},\dots,X_{\lfloor nT\rfloor,n}$, then we can define a random element in $\PP(\XX)$ by
	$$\eta^n := \eta^n((X_{i,n})_{i\le \lfloor nT\rfloor}).$$
\end{itemize}

At this point, it is important to emphasize that the randomness of $\eta^n$ stems only from the random data points and not from the initial configuration of the weights as it is the case for $ \theta^n$ because $\eta^n$ projected at time $0$ is the deterministic distribution $\nu_0^n$. Further, we would like to stress that $\eta^n$ projected at time $t\in[0,T]$ has the same distribution as $\nu_{\lfloor nt\rfloor}^n((X_{i,n})_{i\le \lfloor nt\rfloor})$.

When considering large deviations, we point out that the LDP in Theorem \ref{thm:annealed} for $ \theta^n$ goes beyond large deviations for $\eta^n$ and allows rare occurrences to be caused by a tilted distribution of the initial weights in addition to tilted data. We now present a quenched LDP for $(\eta^n)_n$ in which only the data trajectory measure is tilted.

\begin{theorem}[Quenched LDP for $\eta^n$]
\label{thm:quenched}
Assume that \ref{CON}, \ref{DCOMP}, \ref{UNQ} and \ref{WCOMP'} are satisfied. Let $(\nu_0^n)_n\subseteq \PP(\R^d)$ converge weakly to $\nu$. Then, the family of random probability measures $(\eta^n)_{n \ge 1}$ satisfies the LDP in $\PP(\XX)$ with respect to the weak topology with good rate function $I_\nu$.
\end{theorem}

\section{Proof overview}\label{sec:proof_overview}

\subsection{Proof of Theorem \ref{thm:quenched} -- Quenched LDP}\label{sec:proof_overview_quenched}

The proof of Theorem \ref{thm:quenched} will need three ingredients following the lines of the blueprint for the weak convergence approach in \cite{dupuis}, namely,
\begin{enumerate}
\im[1.] A representation formula before the limit;
\im[2.] Tightness to deduce the existence of a subsequential limit;
\im[3.] Identification of the limit, where we derive the upper and lower bounds.
\end{enumerate}

As in \cite{dupuis}, we consider Laplace principles for the proof of the LDP. In our setting, this means that $(\eta^n)_n$ fulfills the Laplace principle on $\PP(\XX)$ with rate function $I_\nu$ if for every bounded and continuous function $G\colon\PP(\XX)\to\R$ the Laplace upper bound
\begin{equation}\label{eq:laplace_upper}
    \limsup_{n\to\ff}\f1{n'}\log\E\big[\exp\big(-n'G(\eta^n)\big)\big] \le -\inf_{\eta\in\PP(\XX)} \big(G(\eta) + I_\nu(\eta)\big)
\end{equation}
and the Laplace lower bound
\begin{equation}\label{eq:laplace_lower}
    \liminf_{n\to\ff}\f1{n'}\log\E\big[\exp\big(-n'G(\eta^n)\big)\big] \ge -\inf_{\eta\in\PP(\XX)} \big(G(\eta) + I_\nu(\eta)\big)
\end{equation}
hold. As elaborated in \cite[Section 1.2]{dupuis}, a Laplace principle is just an equivalent formulation of the LDP if the potential rate function is good, i.e., if it has compact level sets.

\subsubsection{Representation formula}\label{sec:rep}

The rough idea for the representation in the proof of the quenched LDP is that the original training data is replaced by a tilted variant with a different distribution. To make our presentation self-contained, we repeat here some notation from \cite[Chapter 3]{dupuis}. For $n\in\N$ and a series of probability measures $(\pi_{k,n})_k\subseteq\PP(\R^d)$, let $\bar X_{k, n}$ be distributed according to $\pi_{k,n}$ for each $k\le n T$.
In the case where we plug these random data points into the generalized update equation \eqref{eq:evt2}, we define
$$\bar\eta^n := \eta^n((\bar X_{i,n})_{i\le \lfloor nT\rfloor}).$$
As before, the random data points turn this into a random element of $\PP(\XX)$. Having introduced the notation, we can now restate a modified representation theorem from \cite[Proposition 3.1]{dupuis}. To ease notation, we henceforth set $n' := \lfloor nT\rfloor$.

\begin{proposition}[Representation formula]
\label{thm:rep}
For every bounded measurable function $G\colon \PP(\XX)\to\R$, it holds that 
\begin{equation}\label{eq:rep}
    -\f1{n'}\log\E\big[\exp\big(-n'G(\eta^n) \big)\big] = \inf_{(\pi_{k,n})_k}\E\Big[G(\bar\eta^n) + \f1{n'}\sum_{k\le nT} H(\pi_{k,n} \mid \pi) \Big].
\end{equation}
Here, the infimum is formed over all collections of random probability measures $(\pi_{k,n})_k\subseteq\PP(\R^d)$ such that
\begin{enumerate}
\im the random measure $\pi_{k,n}$ is measurable with respect to the $\s$-algebra $\FF_{k-1,n} = \s((\bar X_{k', n})_{k ' \le k-1})$. Note that we set $\FF_{0, n} := \s(1)$.
\im $\pi_{k,n}$ is the conditional distribution of $\bar X_{k,n}$ given $\FF_{k-1, n}$.
\end{enumerate}
Moreover, for any such collection $(\pi_{k,n})_k$, we have
$$\E\Big[\f1{n'}\sum_{k\le nT} H(\pi_{k,n} \mid \pi) \Big] = \frac1{T} \E[R(\tilde\pi^n)],$$
where
$$\tilde\pi^n({\rm d} t, {\rm d} x) := \sum_{k \le nT}\one\{t \in [(k-1)/n, k/n]\}{\rm d} t\,\pi_{k,n}({\rm d} x)$$
is a random element of $\MM$.
\end{proposition}

\begin{proof}
First, we point out that \eqref{eq:rep} is a direct consequence of  \cite[Proposition 3.1]{dupuis} if we apply it to the bounded function $G':\R^{dn'} \to \R$ given by $G'\big((x_k)_k\big) := G\big(\eta^n((x_k)_k)\big)$. Further, using the chain rule for the relative entropy \cite[Theorem 2.6]{dupuis} we can compute
\begin{align*}
    \E[R(\tilde\pi^n)] &= \E\Big[H(\lambda_T \mid \lambda_T) + \int_{[0,T]} H(\tilde\pi_t^n(\cdot) \mid \pi){\rm d} t \Big] = \E\Big[\sum_{k\le nT} \int_{(k-1)/n}^{k/n} H(\pi_{k,n} \mid \pi){\rm d} t \Big] \\
    &= \E\Big[\frac1{n}\sum_{k\le nT} H(\pi_{k,n} \mid \pi)\Big] = T\E\Big[\frac1{n'}\sum_{k\le nT} H(\pi_{k,n} \mid \pi)\Big],   
\end{align*}
where $\lambda_T$ denotes the Lebesgue measure on $[0,T]$.
\end{proof}

Next, for a bounded measurable function $G\colon\PP(\XX)\to\R$ and every $n\in\N$, we henceforth assume that $(\pi_{k,n})_k$ is a collection of random probability measures from the infimum of Proposition \ref{thm:rep} such that
\begin{equation}\label{eq:bound_entropy}
\sup_{n\in\N} \E\Big[\f1{n'}\sum_{k\le nT} H(\pi_{k,n} \mid \pi)\Big] \le 2\|G\|_\ff + 1.
\end{equation}
This guarantees almost sure absolute continuity with respect to the data distribution.

\subsubsection{Tightness}

The following lemmas deal with the tightness of the sequence of measures $(\tilde\pi^n,\bar\eta^n)_n$ constructed with common sequences of random distributions $((\pi_{k,n})_k)_n$ derived from Proposition \ref{thm:rep}.

\begin{lemma}[Tightness of $\tilde\pi^n$]
\label{lem:tr}
If \eqref{eq:bound_entropy} holds, then the random measures $(\tilde\pi^n)_n$ are tight under $\P$.
\end{lemma}

\begin{lemma}[Tightness of $\bar\eta^n$]
\label{lem:tt}
If \ref{CON}, \ref{DCOMP} and \ref{WCOMP'} are satisfied, then, the random probability measures $(\bar\eta^n)_n$ are tight under $\P$.
\end{lemma}

In contrast to \cite{bufi}, we have to deal with the Skorokhod space and therefore consider the precise dynamics of the training via SGD to prove Lemma \ref{lem:tt}. Also compared to \cite{siri} our proof of tightness of $(\bar\eta^n)_n$ is more involved, as they only have to consider the data trajectory measure $\pi_T$ and not tilted versions.

\subsubsection{Identification -- upper bound}

As a first step, we state that the weak limit of $(\tilde\pi^n,\bar\eta^n)_n$ constructed from the infimum in Proposition \ref{thm:rep} corresponds to a weak solution of \eqref{eq:SDE}.

\begin{lemma}[Weak limit satisfies SDE]\label{lem:weak_convergence_upper_bound}
Assume that $(\tilde\pi^n,\bar \eta^n)$ converges weakly to some $(\tilde\pi,\bar\eta)$ and that
\begin{equation}\label{eq:uniformly_bounded_entropy}
\sup_{n\in\N} \E[R(\tilde\pi^n)] < \ff.
\end{equation}
If \ref{CON}, \ref{DCOMP} and \ref{WCOMP'} are satisfied, then, it holds that almost surely $(\tilde\pi,\bar\eta)\in\PP_\ff^\nu$.
\end{lemma}

We follow the strategy of \cite[Section 6]{kill}, which studies large deviations for an interacting particle system in the setting of continuous functions with respect to the uniform norm topology, where each particle is driven by independent Brownian motions. In contrast, our Skorokhod setting with discrete updates.
\begin{proof}[Proof of Theorem \ref{thm:quenched} -- upper bound]
We start by letting $G$ be a bounded and continuous (with respect to the weak topology) map from $\PP(\XX)$ to $\R$.
Considering an arbitrary $\e>0$, we can find a collection of random probability measures $(\tilde\pi_{k,n})_{k,n}$ from the infimum in Proposition \ref{thm:rep} such that
\begin{align*}
\f1{n'}\log\E\big[\exp\big(-n'G(\eta^n)\big)\big] &\le -\E\Big[G\big(\bar \eta^n\big) + \f1{n'}\sum_{k\le nT} H(\tilde\pi_{k,n} \mid \pi)\Big] + \e = -\E\big[G(\bar \eta^n) + \frac1{T}R(\tilde\pi^n)\big] +\e,
\end{align*}
which in particular implies that $\sup_{n\in\N} \E[ R(\tilde\pi^n)] < \ff$ in accordance with \eqref{eq:bound_entropy}.
In Lemma \ref{lem:tr} and Lemma \ref{lem:tt} we showed tightness of $(\tilde\pi^n, \bar \eta^n)_n$.
Further, the function $U$ given by $U(\tilde\pi^n, \bar \eta^n):= G(\bar \eta^n) + \frac1{T}R(\tilde\pi^n)$ is lower semicontinuous as the sum of a continuous and a lower semicontinuous function. In particular this means that we can choose a subsequence $(\tilde\pi^{n_l},\bar \eta^{n_l})_l$ of $(\tilde\pi^n,\bar \eta^n)_n$ such that $\liminf_{n\to\infty} \E[U(\tilde\pi^n,\bar \eta^n)] = \lim_{l\to\ff} \E[U(\tilde\pi^{n_l},\bar \eta^{n_l})]$. By tightness, the subsequence $(\tilde\pi^{n_l},\bar \eta^{n_l})_l$ contains a further subsequence $(\tilde\pi^{n_{l_z}},\bar \eta^{n_{l_z}})_z$ such that $(\tilde\pi^{n_{l_z}},\bar \eta^{n_{l_z}})$ converges weakly to $(\tilde\pi,\bar\eta)$ for some random element $(\tilde\pi,\bar\eta)\in\MM\times\PP(\XX)$ as $z\to\ff$. Accordingly, we arrive at
$$\liminf_{n\to\infty} \E\big[\big(G(\bar \eta^n) + \frac1{T}R(\tilde\pi^n)\big)\big] = \lim_{l\to\ff} \E\big[\big(G(\bar \eta^{n_l}) + \frac1{T}R(\tilde\pi^{n_l})\big)\big] \ge \E\big[G(\bar \eta) + \frac1{T}R(\tilde\pi)\big].$$
From Lemma \ref{lem:weak_convergence_upper_bound}, we deduce that $(\tilde\pi,\bar\eta)\in\PP_\ff^\nu$ almost surely and thus,
\begin{align*}
\limsup_{n\to\ff} \f1{n'}\log\E\big[\exp\big(-n'G(\eta^n)\big)\big] &\le -\E\big[G(\bar \eta) + \frac1{T} R(\tilde\pi)\big] + \e \le -\E\Big[G(\bar\eta) + \inf_{{(\rho,\eta)}\in\PP_\ff^\nu \colon {\eta} = \bar\eta}\frac1{T} R({\rho})\Big] + \e \\
&\le -\inf_{{\eta}\in\PP(\XX)} \Big(G({\eta}) + \inf_{{(\rho,\tilde\eta)}\in\PP_\ff^\nu \colon {\tilde\eta} = {\eta}} \frac1{T} R({\rho})\Big) + \e.
\end{align*}
Letting $\e\to 0$ concludes the Laplace upper bound.
\end{proof}

\subsubsection{Identification -- lower bound}\label{sec:proof_overview_lower_bound}

For the lower bound we proceed as \cite[Section 7]{kill} but the details of our proof are very distinct from \cite{kill} due to the different settings. We fix a bounded and continuous $G$ mapping from $\PP(\XX)$ to $\R$ and let $(\rho,\bar\eta')\in\PP_\ff^\nu$ be arbitrary.
Next, for $n\in\N$ and $s\in[0,\lfloor nT\rfloor/n]$ define a stochastic kernel by
$$\rho_s^n({\rm d} x) = \sum_{k\le nT} \one\{s\in[(k-1)/n,k/n]\} \int_{(k-1)/n}^{k/n} n \rho_t({\rm d} x){\rm d} t$$

The following lemma asserts that the resulting data distribution converges weakly.

\begin{lemma}[Weak convergence of constructed data distribution]\label{lem:weak_convergence_lower_bound}
The sequence $({\rm d} s\otimes\rho_s^n({\rm d} x))_n$ converges weakly to $\rho$.
\end{lemma}

To prove the Laplace lower bound, we also need to examine the relative entropy of $\rho^n({\rm d} s,{\rm d} x):= {\rm d} s \otimes \rho_s^n({\rm d} x)$ with respect to $\pi_T$.

\begin{lemma}[Relative entropy bounded by weak limit]\label{lem:entropy_bound}
It holds that $\sup_{n\in\N} R(\rho^n) \le R(\rho)$.
\end{lemma}

Next, we prove the lower bound.

\begin{proof}[Proof of Theorem \ref{thm:quenched} -- lower bound]
We can use the measure ${\rm d} s\otimes\rho_s^n({\rm d} x)$ to construct a sequence $(\pi_{k,n})_k$ of probability measures on $\XX$ in the infimum of \eqref{eq:rep}. In particular, for $k\le nT$, we can define
$$\pi_{k,n}({\rm d} x) = \int_{(k-1)/n}^{k/n} n \rho_t^n({\rm d} x){\rm d} t,$$
and let $\bar X_{k,n} \sim \pi_{k,n}$.
With these, we can construct $\bar\eta^n$ as described in Section \ref{sec:rep}, where we set $\nu_0^n = \nu$.
We also know from Lemma \ref{lem:tr}, which becomes applicable because of the entropy bound in Lemma \ref{lem:entropy_bound}, and Lemma \ref{lem:tt} that $(\rho^n,\bar\eta^n)$ is tight and every subsequence has a subsubsequence that converges to $(\rho,\bar\eta)$ for some $\bar\eta$. By Lemma \ref{lem:weak_convergence_upper_bound} this limit has to be in $\PP_\ff^\nu$. Now, assumption \ref{UNQ} guarantees that all these subsequential limits are equal and thus, weak convergence of $(\rho^n,\bar\eta^n)$ to $(\rho,\bar\eta')$ follows.
After pointing out that $\frac1{T} R(\rho^n) = \f1{n'}\sum_{k\le nT} H(\pi_{k,n} \mid \pi)$, this lets us conclude that

\begin{align*}
\liminf_{n\to\ff} \f1{n'}\log\E\big[\exp\big(-n'G(\eta^n)\big)\big] &\ge \liminf_{n\to\ff} -\Big(\E\big[G(\bar \eta^n)\big] + \f1{n'}\sum_{k\le nT} H(\pi_{k,n} \mid \pi)\Big) \\
&\ge -\E\big[G(\bar\eta')\big] - \frac1{T}R(\rho) = -G(\bar\eta') - \frac1{T} R(\rho),
\end{align*}
where the first inequality follows from Proposition \ref{thm:rep} and the second one from weak convergence of $(\rho^n,\bar\eta^n)$ to $(\rho,\bar\eta')$ and Lemma \ref{lem:entropy_bound}.
Invoking that $(\rho,\bar\eta')\in\PP_\ff^\nu$ was arbitrary completes the argument.
\end{proof}

\subsubsection{Identification -- good rate function}

This section is devoted to presenting the result that $I_\nu$ is indeed a good rate function, i.e., $I_\nu$ has compact level sets. We even prove a more sophisticated version of that because of its usefulness in the proof of Theorem \ref{thm:annealed}, when transitioning from the quenched LDP to the annealed LDP.

\begin{lemma}[Compact sublevel sets]\label{lem:compact_level_sets}
Assume that \ref{CON}, \ref{DCOMP} and \ref{WCOMP'} are satisfied. Then, for any $M<\infty$ the set $\{(\nu',{\eta})\in\PP(\{(c,w)\in\R^d: \|(c,w)\|\le C_\nu\})\times\PP(\XX)\colon I_{\nu'}({\eta})\le M\}$ is compact with respect to the weak topology.
\end{lemma}

The following corollary follows directly follows from Lemma \ref{lem:compact_level_sets}.

\begin{corollary}\label{cor:compact_level_sets}
Assume that \ref{CON}, \ref{DCOMP} and \ref{WCOMP'} are satisfied. Then, for any $M<\infty$ the set $\{{\eta}\in\PP(\XX)\colon I_\nu({\eta})\le M\}$ is compact with respect to the weak topology.
\end{corollary}

By Corollary \ref{cor:compact_level_sets} and since $I_\nu$ is a nonnegative function it is a good rate function. Together with the upper and lower bounds from the previous sections, this proves a Laplace principle with good rate function for $\eta^n$ and thus, concludes the proof of Theorem \ref{thm:quenched}.

\subsection{Proof of Theorem \ref{thm:annealed} -- Annealed LDP}\label{sec:proof_overview_annealed}

In this section, we show how to derive the annealed LDP from the quenched LDP.
The following lemma asserts that the construction of empirical weight trajectories is indeed a special case of weight trajectories that start from an arbitrary distribution.

\begin{lemma}[Weak convergence of initial weight measures]\label{lem:empirical}
Let $\nu_0^n := \frac1n\sum_{i=1}^n \delta_{\theta_0^{i,n}}$ be the random element given after \eqref{eq_g}. Then, $\eta^n = \theta^n$ in distribution and $\nu_0^n$ converges weakly to $\nu$.
\end{lemma}

\begin{proof}
Fix data points $x_{1,n},\dots,x_{\lfloor nT\rfloor,n}\in\R^d$. By construction of the trajectory distribution in \eqref{eq:c2} and \eqref{eq:evt2} we see that $\ms{Law}_{\nu_0^n}(\tilde\theta_{1/n}^n(X_{1,n})) = \frac1n\sum_{i=1}^n \delta_{\theta_{1/n}^{i,n}}$ in distribution (with respect to $\P$). By induction, we recover that $\eta^n = \theta^n$ in distribution. That $\nu_0^n$ converges weakly to $\nu$, is implied by the Glivenko-Cantelli theorem \cite[Lemma 3.2]{dupuis}, which states that the empirical distribution function of independent and identically distributed random variables converges almost surely uniformly to the distribution function of their common distribution. Thus, in particular the empirical measure $\frac1n\sum_{i=0}^n \delta_{\theta_0^{i,n}}$ converges weakly to $\nu$.
\end{proof}

Under certain conditions annealed large deviations follow from the quenched large deviations. \cite[Theorem 2.3]{annealedLDP} states a set of requirements that enables this transition and which has been successfully applied in the past for instance in \cite{Ganesh}. Those conditions adapted to our setting are the following:
\begin{enumerate}
    \item $\PP(B_{C_\nu}(0))$ and $\PP(\XX)$ are Polish spaces.
    \item $(\frac1n \sum_{i=1}^n \delta_{\theta_0^{i,n}})_n$ satisfies the LDP in $\PP(B_{C_\nu}(0))$ with respect to the weak topology with good rate function $H(\,\cdot\mid\nu)$.
    \item For each $\nu\in\PP(B_{C_\nu}(0))$, the function $I_\nu$ is a good rate function, i.e., it is nonnegative, lower semicontinuous and has compact level sets.
    \item For any sequence $(\nu_0^n)_n\subseteq\PP(B_{C_\nu}(0))$ converging weakly to $\nu$, the random probability measure $\eta^n$ satisfies the LDP in $\PP(\XX)$ with respect to the weak topology with rate function $I_\nu$.
    \item The function $(\rho,{\eta}) \mapsto I_\rho({\eta})$ is lower semicontinuous.
\end{enumerate}
Verifying these is the next step.

\begin{proof}[Proof of Theorem \ref{thm:annealed}]
\begin{enumerate}
    \item $\XX$ equipped with the Skorokhod topology is a Polish space as a consequence of \cite[Theorem 5.6]{ethierkurtz1986}, which states that a Skorokhod space is complete and separable if the underlying function domain is complete and separable. By \cite[Lemma 4.5]{randMeas}, this turns $\PP(\XX)$ equipped with the weak topology into a Polish space. The latter also implies that $\PP(B_{C_\nu}(0))$ as a closed subset of $\PP(\R^d)$ is Polish when equipped with the weak topology.
    \item This is Sanov's theorem, see for example \cite[Theorem 6.2.10]{dz98}, which deals with the LDP for empirical laws on Polish spaces.
    \item By Corollary \ref{cor:compact_level_sets} and since $I_\nu$ is a nonnegative function for all $\nu$, the function $I_\nu$ is a good rate function.
    \item This is the quenched LDP from Theorem \ref{thm:quenched}.
    \item Finally, by Lemma \ref{lem:compact_level_sets} the function $(\rho,{\eta}) \mapsto I_\rho({\eta})$ is lower semicontinuous.
\end{enumerate}

These observations together with Lemma \ref{lem:empirical} let us apply \cite[Theorem 2.3]{annealedLDP}, which implies that $\theta^n$ satisfies the LDP with the stated rate function.
\end{proof}

\section{Proof of tightness}
\label{sec:tight}

Now, we show tightness of the sequences measures $(\tilde\pi^n)_n$ and $(\bar\eta^n)_n$. That $(\tilde\pi^n)_n$ is tight is established via a tightness function argument.

\begin{proof}[Proof of Lemma \ref{lem:tr}]
Note that $(\tilde\pi^n)_n$ is a sequence of random elements in $\MM$. By \cite[Theorem 2.11]{dupuis}, tightness of the collection of intensity measures $(\E[\tilde\pi^n])_n$ implies that $(\tilde\pi^n)_n$ is tight. Next, \cite[Lemma 2.4 (c)]{dupuis} and \cite[Lemma 2.10]{dupuis} assert that the relative entropy $R$ is a tightness function on $\MM$, i.e., it is nonnegative and has precompact level sets. Using \cite[Lemma 2.9]{dupuis}, it now suffices to show that $\sup_{n \in \N}\int_\MM R(\rho)\E[\tilde\pi^n]({\rm d}\rho) < \ff$. Note that we can decompose each realization $\rho$ of $\tilde\pi^n$ into a product measure of the Lebesgue measure $\lambda_T$ on $[0,T]$ and a stochastic kernel on $\R^d$ given $[0,T]$ that we denote by $\rho_t$. Now, using the chain rule for the relative entropy \cite[Theorem 2.6]{dupuis} we arrive at
\begin{align*}
	\int_\MM R(\rho)\E[\tilde\pi^n]({\rm d}\rho) &= \E\Big[\int_\MM R(\rho) \tilde\pi^n({\rm d}\rho)\Big] = \E\Big[\int_\MM \Big(H(\lambda_T \mid \lambda_T) + \int_{[0,T]} H(\rho_t \mid \pi){\rm d} t\Big) \tilde\pi^n({\rm d}\rho)\Big] \\
	&= \E\Big[ \int_\MM \sum_{k\le nT} \int_{(k-1)/n}^{k/n} H(\pi_{k,n} \mid \pi){\rm d} t \tilde\pi^n({\rm d}\rho)\Big] = \E\Big[\frac1{n}\sum_{k\le nT} H(\pi_{k,n} \mid \pi)\Big] \\
	&= \frac{n'}{n}\E\Big[\frac1{n'}\sum_{k\le nT} H(\pi_{k,n} \mid \pi)\Big] \le T(2\|G\|_\ff + 1).
\end{align*}
\end{proof}

It remains to prove Lemma \ref{lem:tt}. Here, we aim to show that with high probability $(\bar\eta^n)_n$ can be reduced to a measure on a compact space by using the characterization of compactness in the Skorokhod space via continuity modulus given in \cite[Theorem 12.3]{billingsley}.

For this, we define possible weight trajectory with fixed potential initial values $\omega_0\in\supp(\nu_0^n)$, by recursively setting
$$
	\hat\theta_{k/n}^n((x_{i,n})_{i\le k}) :=\ms{SGD}(\hat\theta_{(k-1)/n}^n((x_{i,n})_{i\le k-1}), x_{k,n}, \nu_{k-1}^n((x_{i,n})_{i\le k-1})),
$$
where $\hat\theta_0^n := \omega_0$. This also gives a $\XX$-valued random variable $(\hat\theta_t^n)_t$ if we plug in random data points, i.e., by setting
\begin{equation}\label{eq:th_def}
	\hat\theta_t^n := \hat\theta_t^n((\bar X_{i,n})_{i\le \lfloor tn\rfloor}) := \hat\theta_{\lfloor tn\rfloor/n}^n((\bar X_{i,n})_{i\le \lfloor tn\rfloor})
\end{equation}
for $t\in[0,T]$. We omit the dependence on $\omega_0$ in the notation.

As an essential preliminary step, we derive growth bounds on $\hat\theta^{n}_t$ depending on the tilted data points. These are derived by carefully analyzing the update equation \eqref{eq:c2}. To simplify the notation, we set $\bar Y_n^{*,m} := \frac1n\sum_{i=1}^{\lfloor nT\rfloor} |\bar Y_{i,n}|^m$ and $\bar Z_n^{*,m} := \frac1n\sum_{i=1}^{\lfloor nT\rfloor} \|(1,\bar Z_{i,n})\|^m$ for $m\in\N$.

\begin{lemma}[Growth bound on $(\hat\theta_t^n)_t$]
\label{lem:th_gr}
Assume that \ref{CON} and \ref{WCOMP'} hold. Then,
\begin{enumerate}[leftmargin=*]
\im[(i)]  $\sup_{\omega_0\in B_{C_\nu}(0)} \frac1n \sum_{i=1}^{\lfloor nT\rfloor} \|A(\bar X_{k,n}, \hat\theta_{(k-1)/n}^n; \bar\eta_{(k-1)/n}^n)\| \le C_\ms{SGD} (T^2+1)\big(1 + \bar Y_n^{*,4} + \bar Z_n^{*,2}\big)$ for some $C_\ms{SGD}>1$,
\im[(ii)]  $\sup_{\omega_0\in B_{C_\nu}(0)} \|(\hat\theta_t^n)_t\|\le C_\nu + C_\ms{SGD} (T^2+1)\big(1 + \bar Y_n^{*,4} + \bar Z_n^{*,2}\big)$
\end{enumerate}
and if additionally \ref{DCOMP} holds, then
\begin{enumerate}[leftmargin=*]
\im[(iii)]  $\lim_{M\to\ff} \sup_{n\in\N} \E[\frac1n\sum_{k=1}^{\lfloor nT\rfloor} |\bar Y_{k,n}|^{{4}} \one_{\{|\bar Y_{k,n}| \ge M\}}] = 0$ and \\
$\lim_{M\to\ff} \sup_{n\in\N} \E[\frac1n\sum_{k=1}^{\lfloor nT\rfloor} \|(1,\bar Z_{k,n})\|^{{2}} \one_{\{\|(1,\bar Z_{k,n})\|\ge M\}}] = 0$. \\
In particular, $$\lim_{M\to\ff} \sup_{n\in\N} \sup_{\omega_0\in B_{C_\nu}(0)} \E\bigg[\frac1n\sum_{k=1}^{\lfloor nT\rfloor} \|A(\bar X_{k,n}, \hat\theta_{(k-1)/n}^n; \bar\eta_{(k-1)/n}^n)\| \one_{\{\|\bar X_{k,n}\|\ge M\}}\bigg] = 0.$$
\end{enumerate}
\end{lemma}

\begin{proof}
We start by defining $(c_{k/n}^n, w_{k/n}^n) := \hat\theta_{k/n}^n$ for each $k\le nT$. Since $\s$ as well as $\s'$ are bounded by \ref{CON}, it follows that
\begin{align}\label{eq:bd_loss}
\begin{split}
|g_{k,n}| &:= |g\big((\bar Z_{k,n}, \bar Y_{k,n}), \bar\eta_{(k-1)/n}^n\big)| \le |\bar Y_{k,n}| + \int_{\R^d} |c \s(\bar Z_{k,n}^\top w)| \bar\eta_{(k-1)/n}^n {\rm d}(c,w) \\
&\le |\bar Y_{k,n}| + C_\s \int_{\R^d} |c| \bar\eta_{(k-1)/n}^n {\rm d}(c,w).
\end{split}
\end{align}
Now, by the definition of SGD in \eqref{eq:c2} and \eqref{eq:grad}, and with \eqref{eq:bd_loss} we can bound
\begin{align*}
\int_{\R^d} |c| \bar\eta_{k/n}^n {\rm d}(c,w) &\le \int_{\R^d} |c| + \frac1n |g_{k,n} \s(\bar Z_{k,n}^\top w)| \bar\eta_{(k-1)/n}^n {\rm d}(c,w) \le \int_{\R^d} |c| \bar\eta_{(k-1)/n}^n {\rm d}(c,w) + \frac1n C_\s |g_{k,n}| \\
&\le C_\s |\bar Y_{k,n}| / n + (1 + C_\s^2 /n) \int_{\R^d} |c| \bar\eta_{(k-1)/n}^n {\rm d}(c,w).
\end{align*}
By induction, $a_k \le C' b_k + C a_{k-1}$ for $C,C'>0$ leads to $a_k \le C^k a_0 + \sum_{i=1}^k C^{k-i-1} C' b_i $ and thus,
\begin{align}\label{eq:bd_average_c}
\begin{split}
\int_{\R^d} |c| \bar\eta_{k/n}^n {\rm d}(c,w) &\le (1 + C_\s^2/n)^k \int_{\R^d} |c| \nu_0^n {\rm d}(c,w) + \sum_{i=1}^k (1 + C_\s^2/n)^{k-i-1} C_\s |\bar Y_{i,n}| /n  \\
&\le C_\nu e^{C_\s^2 T} + C_\s e^{C_\s^2 T} \bar Y_n^{*,1},
\end{split}
\end{align}
where we used that $(1+C/n)^k \le (1+C/n)^{nT} \le e^{CT}$ for $C>0$ and that $\nu_0^n$ has compact support by \ref{WCOMP'}. Further, note that combining \eqref{eq:bd_loss} and \eqref{eq:bd_average_c} with the update equation \eqref{eq:c2} yields
\begin{align*}
|c_{k/n}^n| &\le |c_{(k-1)/n}^n| + \frac1n |g_{k,n} \s(\bar Z_{k,n}^\top w_{(k-1)/n}^n)| \le  |c_{(k-1)/n}^n| + \frac1n \Big(|\bar Y_{k,n}| + C_\s \int_{\R^d} |c| \bar\eta_{(k-1)/n}^n {\rm d}(c,w)\Big) C_\s \\
&\le  |c_{(k-1)/n}^n| + \frac1n \big(|\bar Y_{k,n}| + C_\s \big(C_\nu e^{C_\s^2 T} + C_\s e^{C_\s^2 T} \bar Y_n^{*,1} \big)\big) C_\s.
\end{align*}
Thus, we arrive at
\begin{equation}\label{eq:bd_c}
|c_{k/n}^n| \le |c_0^n| + C_\s \bar Y_n^{*,1} + T C_\s^2 C_\nu e^{C_\s^2 T} + T C_\s^3 e^{C_\s^2 T} \bar Y_n^{*,1} \le \underbrace{2(1+T) C_\s^3 e^{C_\s^2 T}}_{=:\bar C} \big(C_\nu + \bar Y_n^{*,1}\big).
\end{equation}
Next, we bound the weight evolution in an SGD step for a single neuron by
\begin{align}\label{eq:update_bound}
\begin{split}
\|A(\bar X_{k,n}, \hat\theta_{(k-1)/n}^n; \bar\eta_{(k-1)/n}^n)\| &= \big\|\big(g_{k,n} \s(\bar Z_{k,n}^\top w_{(k-1)/n}^n),g_{k,n} c_{(k-1)/n}^n \s'(\bar Z_{k,n}^\top w_{(k-1)/n}^n) \bar Z_{k,n} \big)\big\| \\
&\le |g_{k,n}| C_\s \big\|\big(1,c_{(k-1)/n}^n \bar Z_{k,n} \big)\big\| \\
&\le \big(|\bar Y_{k,n}| + C_\s^2 C_\nu e^{C_\s^2 T} + C_\s^3 e^{C_\s^2 T} \bar Y_n^{*,1} \big) \bar C \big(C_\nu + \bar Y_n^{*,1}\big) \|(1,\bar Z_{k,n})\|,
\end{split}
\end{align}
From this, it follows that
\begin{align*}
&\frac1n \sum_{i=1}^{\lfloor nT\rfloor} \|A(\bar X_{i,n}, \hat\theta_{(i-1)/n}^n; \bar\eta_{(i-1)/n}^n)\| \\
&\le \frac1n\sum_{i=1}^{\lfloor nT\rfloor} \big(|\bar Y_{i,n}| + C_\s^2 C_\nu e^{C_\s^2 T} + C_\s^3 e^{C_\s^2 T} \bar Y_n^{*,1} \big) \bar C \big(C_\nu + \bar Y_n^{*,1}\big) \|(1,\bar Z_{i,n})\| \\
&\le \bar C C_\nu \frac1n\sum_{i=1}^{\lfloor nT\rfloor} |\bar Y_{i,n}| \|(1,\bar Z_{i,n})\| + \bar C \bar Y_n^{*,1} \big(\frac1n\sum_{i=1}^{\lfloor nT\rfloor} |\bar Y_i| \|(1,\bar Z_{i,n})\|\Big) + \bar C  C_\s^2 C_\nu^2 e^{C_\s^2 T} \bar Z_n^{*,1} \\
&\quad+ 2 \bar C C_\s^3 C_\nu e^{C_\s^2 T} \bar Y_n^{*,1} \bar Z_n^{*,1} + \bar C C_\s^3 e^{C_\s^2 T} \big(\bar Y_n^{*,1}\big)^2 \bar Z_n^{*,1}.
\end{align*}
Now, we use Cauchy-Schwarz inequality to get multiple bounds. First, we see that
\begin{align*}
	\frac1n\sum_{i=1}^{\lfloor nT\rfloor} |\bar Y_{i,n}| \|(1,\bar Z_{i,n})\| \le \frac1n \sqrt{\sum_{i=1}^{\lfloor nT\rfloor} |\bar Y_{i,n}|^2} \sqrt{\sum_{i=1}^{\lfloor nT\rfloor} \|(1,\bar Z_{i,n})\|^2} \le \bar Y_n^{*,2} + \bar Z_n^{*,2}.
\end{align*}
Second, it holds that
\begin{align*}
	\bar Y_n^{*,1} \bar Z_n^{*,1} &\le \frac1{n^2} \Big(\sum_{i=1}^{\lfloor nT\rfloor} |\bar Y_{i,n}|\Big)^2 + \frac1{n^2} \Big(\sum_{i=1}^{\lfloor nT\rfloor} \|(1,\bar Z_{i,n})\|\Big)^2 \le T \bar Y_n^{*,2} + T \bar Z_n^{*,2} \\
	&\le T \sqrt{T Y_n^{*,4}} + T \bar Z_n^{*,2} \le T + T^2 \bar Y_n^{*,4} + T \bar Z_n^{*,2}
\end{align*}
and third, we have that
\begin{align*}
	(\bar Y_n^{*,1})^2 \bar Z_n^{*,1} &\le T \bar Y_n^{*,2} \bar Z_n^{*,1} \le T \bar Y_n^{*,4} + T \bar Z_n^{*,2}.
\end{align*}
Combining these bounds yields for some $C_\ms{SGD}>0$
$$
\frac1n \sum_{i=1}^{\lfloor nT\rfloor} \|A(\bar X_{i,n}, \hat\theta_{(i-1)/n}^n; \bar\eta_{(i-1)/n}^n)\| \le C_\ms{SGD}(T^2+1)\big(1 + \bar Y_n^{*,4} + \bar Z_n^{*,2}\big).
$$
To show part (ii), for $t\in[0,1/n)$ it holds that $\|\hat\theta_t^n\| \le C_\nu$ and for $t\in[1/n, T]$ choose $k\le nT$ such that $k/n \le t < (k+1)/n$. Then, it follows that $\hat\theta_t^n = \hat\theta_{k/n}^n$ and we can compute
\begin{align*}
\|\hat\theta_t^n\| &= \|\hat\theta_{k/n}^n\| = \|\hat\theta_{(k-1)/n}^n + \frac1n A(\bar X_{k,n}, \hat\theta_{(k-1)/n}^n; \bar\eta_{(k-1)/n}^n)\| \le\|\hat\theta_0^n\| + \frac1n\sum_{i=1}^{\lfloor nT\rfloor} \|A(\bar X_{i,n}, \hat\theta_{(i-1)/n}^n; \bar\eta_{(i-1)/n}^n)\| \\
&\le C_\nu + C_\ms{SGD} (T^2+1)\big(1 + \bar Y_n^{*,4} + \bar Z_n^{*,2}\big).
\end{align*}
Finally, for part (iii), we aim to use the following inequality
$$a b \le e^{a h} + \frac1{h} (b\log b - b + 1)$$
for $a,b\ge 0$ and $h\ge 1$ from \cite[Equation (2.9)]{dupuis}. With this and Cauchy-Schwarz inequality we proceed as in \cite[Proof of Lemma 3.9]{dupuis} and compute for any $M,h>0$
\begin{align*}
	&\E\bigg[\frac1n\sum_{k=1}^{\lfloor nT\rfloor} |\bar Y_{k,n}|^{{4}} \one\{|\bar Y_k| \ge M\}\bigg] = \E\bigg[\frac1n\sum_{k=1}^{\lfloor nT\rfloor} \E\big[ |\bar Y_{k,n}|^{{4}} \one\{ |\bar Y_{k,n}| \ge M\} \bigm\vert \FF_{k-1,n}\big] \bigg] \\
	&= \frac1n \E\bigg[\sum_{k=1}^{\lfloor nT\rfloor} \int_{\R^d} |y|^{{4}} \one\{ |y| \ge M\} \pi_{k,n}({\rm d} (z,y)) \bigg] = \frac1n \sum_{k=1}^{\lfloor nT\rfloor} \E\bigg[\int_{\R^d} |y|^{{4}} \one\{ |y| \ge M\} \frac{{\rm d}\pi_{k,n}}{{\rm d}\pi}(z,y) \pi({\rm d} (z,y))\bigg] \\
	&\le \frac1n \sum_{k=1}^{\lfloor nT\rfloor} \int_{\R^d} e^{h |y|^{{4}}}\one\{ |y| \ge M\} \pi({\rm d} (z,y)) + \frac1{h} \E[H(\pi_{k,n}\mid \pi)] \\
	&\le \sqrt{\E[e^{2 h |Y|^{{4}}}]} \sqrt{\P(|Y| \ge M)} + \frac1n \sum_{k=1}^{\lfloor nT\rfloor} \frac1{h} \E[H(\pi_{k,n}\mid \pi)] \le C_h \sqrt{\P(|Y| \ge M)} + \frac1{h} T(2\|G\|_\ff + 1)
\end{align*}
for some $C_h>0$ by \ref{DCOMP} and with \eqref{eq:bound_entropy} used in the last inequality.
When letting $M\to\ff$ first and then $h\to\ff$, this implies that
$$\lim_{M\to\ff} \sup_{n\in\N} \E\bigg[\frac1n\sum_{k=1}^{\lfloor nT\rfloor} |\bar Y_{k,n}|^{{4}} \one\{|\bar Y_{k,n}| \ge M\}\bigg] = 0$$
A similar calculation leads to
$$\lim_{M\to\ff} \sup_{n\in\N} \E\bigg[\frac1n\sum_{k=1}^{\lfloor nT\rfloor} \|(1,\bar Z_{k,n})\|^{{2}} \one_{\{\|(1,\bar Z_{k,n})\|\ge M\}}\bigg] = 0.$$
Using (i), we can conclude that
$$\lim_{M\to\ff} \sup_{n\in\N} \sup_{\omega_0\in B_{C_\nu}(0)} \E\bigg[\frac1n\sum_{k=1}^{\lfloor nT\rfloor} \|A(\bar X_{k,n}, \hat\theta_{(k-1)/n}^n; \bar\eta_{(k-1)/n}^n)\| \one_{\{\|\bar X_{k,n}\|\ge M\}}\bigg] = 0.$$
\end{proof}

We are now ready to prove tightness of the sequence $(\bar\eta^n)_n$. Note that the proof of Lemma \ref{lem:tt} will also yield the following result that is useful later.

\begin{lemma}[Expected maximal jump bound]
\label{lem:max_bd}
Assume that \ref{CON}, \ref{DCOMP} and \ref{WCOMP'} are satisfied. Then, it holds that
 $$\limsup_{n\to\ff} \sup_{\omega_0\in B_{C_\nu}(0)}\E\bigg[\max_{i\le nT} \frac1n \|A(\bar X_{i,n}, \hat\theta_{(i-1)/n}^n; \bar\eta_{(i-1)/n}^n)\| \bigg] = 0.$$
\end{lemma}

\begin{proof}[Proof of Lemma \ref{lem:tt} and Lemma \ref{lem:max_bd}]
We aim to show tightness by using a compactness characterization for sets on a Skorokhod space from \cite[Section 12]{billingsley} based on bounds for the continuity modulus for $(\theta_t)_t\in\XX$. First, for $\delta>0$ we define
$$w'_{(\theta_t)_t}(\delta) := \inf_{\{t_i\}_{0\le i\le l}} \max_{1\le i\le l} \sup_{s,t\in[t_{i-1},t_i)} \|\theta_t - \theta_s\|,$$
where the infimum is taken over all time points $0=t_0 < t_1 < \ldots < t_l = T$ that form a partition of $[0,T]$ of arbitrary size $l\in\N$ but satisfy $t_i - t_{i-1} > \delta$.
We aim to show that for all $\e>0$
\begin{equation}\label{eq:tightness_crit}
	\lim_{\delta\to0}\sup_{n\to\infty} \sup_{\omega_0\in B_{C_\nu}(0)}\P(w'_{(\hat\theta_t^n)_t}(\delta) > \e) = 0.
\end{equation}
To achieve this, first note that if $\delta<1/n$, then $w'_{(\hat\theta_t^n)_t}(\delta) = 0$ for all $\omega_0\in B_{C_\nu}(0)$, since we can choose the time points to be the points at which there is a jump in the sequence $(\hat\theta_t^n)_t$, i.e., $t_i = i/n$ for $i\le nT$.
In the case that $\delta\ge 1/n$ we apply Markov's inequality and assert that for some constant $C'>0$
\begin{align}\label{eq:compact_assertion}
\begin{split}
	\P(w'_{(\hat\theta_t^n)_t}(\delta) \ge \e) &\le \E\bigg[\max_{\substack{0\le k_1<k_2\le nT, \\ k_2 - k_1\in (\delta n,2\delta n + 1]}} \frac1n \sum_{i=k_1}^{k_2} \|A(\bar X_{i,n}, \hat\theta_{(i-1)/n}^n; \bar\eta_{(i-1)/n}^n)\| \bigg] / \e \\
	&\le \kappa/\e + (2\delta + 1/n) C' (1+M) (M + C')/\e \le \kappa/\e + 3\delta C' (1+M) (M + C')/\e
\end{split}
\end{align}
for all $\kappa>0$ if $M := M(\kappa)$ is large enough. Since the bound does not depend on $\omega_0$, letting $\delta\to0$ and $\kappa\to0$ in this particular order yields \eqref{eq:tightness_crit}. We point out that this assertion also proves Lemma \ref{lem:max_bd} because
\begin{align*}
&\E\bigg[\max_{i\le nT} \frac1n \|A(\bar X_{i,n}, \hat\theta_{(i-1)/n}^n; \bar\eta_{(i-1)/n}^n)\| \bigg] \le \E\bigg[\max_{\substack{0\le k_1<k_2\le nT, \\ k_2 - k_1\in (\delta n,2\delta n + 1]}} \frac1n \sum_{i=k_1}^{k_2} \|A(\bar X_{i,n}, \hat\theta_{(i-1)/n}^n; \bar\eta_{(i-1)/n}^n)\| \bigg].
\end{align*}
To show the assertion, note that by Lemma \ref{lem:th_gr} (iii) and with the Cauchy-Schwarz inequality, for each $\kappa>0$, we can choose $M$ large enough such that for every $n\in\N$ and $\omega_0\in B_{C_\nu}(0)$
\begin{align*}
	&\E\bigg[\max_{\substack{0\le k_1<k_2\le nT, \\ k_2 - k_1\in (\delta n,2\delta n + 1]}} \frac1n \sum_{i=k_1}^{k_2} \|A(\bar X_{i,n}, \hat\theta_{(i-1)/n}^n; \bar\eta_{(i-1)/n}^n)\| \one\{ \|\bar X_{i,n}\| \ge M\}\bigg] \\
	&\le \E\bigg[\frac1n \sum_{i=1}^{\lfloor nT\rfloor} \|A(\bar X_{i,n}, \hat\theta_{(i-1)/n}^n; \bar\eta_{(i-1)/n}^n)\| \one\{ \|\bar X_{i,n}\| \ge M\}\bigg] \le \kappa.
\end{align*}
Further, \eqref{eq:update_bound} from the proof of Lemma \ref{lem:th_gr} implies that for $M>0$
\begin{align*}
	&\E\Big[\max_{\substack{0\le k_1<k_2\le nT, \\ k_2 - k_1\in (\delta n,2\delta n + 1]}} \frac1n \sum_{i=k_1}^{k_2} \|A(\bar X_{i,n}, \hat\theta_{(i-1)/n}^n; \bar\eta_{(i-1)/n}^n)\| \one\{\|\bar X_{i,n}\| \le M\} \Big] \\
	&\le \E\Big[\max_{\substack{0\le k_1<k_2\le nT, \\ k_2 - k_1\in (\delta n,2\delta n + 1]}} \frac1n \sum_{i=k_1}^{k_2}\big(|\bar Y_{i,n}| + C_\s^2 C_\nu e^{C_\s^2 T} + C_\s^3 e^{C_\s^2 T} \bar Y_n^{*,1} \big) 2(1+T) C_\s^3 e^{C_\s^2 T} \\
	&\qquad\qquad\qquad\qquad\qquad\qquad\quad\cdot \big(C_\nu + \bar Y_n^{*,1}\big) \|(1,\bar Z_{i,n})\| \one\{\|\bar X_{i,n}\| \le M\}\Big] \\
	&\le \E\Big[\max_{\substack{0\le k_1<k_2\le nT, \\ k_2 - k_1\in (\delta n,2\delta n + 1]}} \frac1n \sum_{i=k_1}^{k_2} \big(M + C_\s^2 C_\nu e^{C_\s^2 T} + C_\s^3 e^{C_\s^2 T} \bar Y_n^{*,1} \big) 2(1+T) C_\s^3 e^{C_\s^2 T} \cdot \big(C_\nu + \bar Y_n^{*,1}\big) (1+M)\Big] \\
	&\le (2\delta + 1/n) C' (1+M) \big(1 + M + \E\big[\bar Y_n^{*,1} \big] + \E\big[T\bar Y_n^{*,2} \big]\big) \\
	&\le (2\delta + 1/n) C' (1+M) (M + C'),
\end{align*}
for some constant $C'>0$, where we used the Cauchy-Schwarz inequality and that the supremum over all $n$ of the two occurring expectations is finite by Lemma \ref{lem:th_gr} (iii). Thus, we have established \eqref{eq:compact_assertion}.

Further, for every $\gamma>0$ we can choose $r>0$ large enough such that
\begin{align*}
	\sup_{\omega_0\in B_{C_\nu}(0)} \P(\|(\hat\theta_t^n)_t\|_\ff > r) &\le \sup_{\omega_0\in B_{C_\nu}(0)} \P\bigg(C_\nu + \frac1n \sum_{i=1}^{\lfloor nT\rfloor} \|A(\bar X_{k,n}, \hat\theta_{(k-1)/n}^n; \bar\eta_{(k-1)/n}^n)\| > r\bigg) \\
	&\le \sup_{n\in\N} \sup_{\omega_0\in B_{C_\nu}(0)} \E\bigg[C_\nu + \frac1n \sum_{i=1}^{\lfloor nT\rfloor} \|A(\bar X_{k,n}, \hat\theta_{(k-1)/n}^n; \bar\eta_{(k-1)/n}^n)\|\bigg] / r \le \gamma,
\end{align*}
since by Lemma \ref{lem:th_gr} (iii) the expectation is bounded uniformly in $n$.

Next, we construct compact sets that include the trajectories $(\hat\theta_t^n)_t$ with high probability and do not depend on $\omega_0$ following the lines of \cite[Theorem 7.2 and Theorem 12.3]{billingsley}, the Arzelà-Ascoli theorem. Let $B_r := \{(\theta_t)_t\in\XX\co\|(\theta_t)_t\|_\ff \le r\}$. We know that for $\gamma>0$ we can choose $r$ large enough such that $\P((\hat\theta_t^n)_t \in B_r) \ge 1-\gamma$ for all $n\in\N$ and all $\omega_0\in B_{C_\nu}(0)$. Further, choose $\delta_m>0$ such that if $B'_m := \{(\theta_t)_t\in\XX\co w'_{(\theta_t)_t}(\delta_m) \le 1/m\}$, then $\P((\hat\theta_t^n)_t \in B'_m) \ge 1 - \gamma/2^m$ for all $n\in\N$ and all $\omega_0\in B_{C_\nu}(0)$. For $\bar K$ defined as the closure of $B_r\cap\bigcap_m B'_m$ it holds that
$$\P((\hat\theta_t^n)_t \in \bar K) \ge 1 - \P((\hat\theta_t^n)_t \not\in B_r) - \sum_{m=1}^\ff \P((\hat\theta_t^n)_t \not\in B'_m) \ge 1 - 2\gamma.$$
By \cite[Theorem 12.3]{billingsley}, $\bar K$ is compact. This means $\PP(\bar K)$, the set of probability measures on $\bar K$ is tight. We may invoke Prokhorov's theorem (\cite[Theorem 23.2]{kallenberg}) to deduce that $\PP(\bar K)$ is compact with respect to the weak topology. Now, we can utilize the similar construction of $(\hat\theta_t^n)_t$ and $\eta^n$ and the fact that $\bar K$ does not depend on $\omega_0$ to conclude that $\P(\bar\eta^n \not\in \PP(\bar K)) \le \sup_{\omega_0\in B_{C_\nu}(0)} \P((\hat\theta_t^n)_t \not\in \bar K) \le 2\gamma$ for all $n\in\N$. Thus, we have shown tightness of $(\bar\eta^n)_n$.
\end{proof}

\section{Proofs for Theorem \ref{thm:quenched} -- quenched LDP}
\label{sec:ident}

\subsection{Proofs -- upper bound}\label{sec:laplace_upper_bound}
It remains to prove Lemma \ref{lem:weak_convergence_upper_bound}. {The crux of this section is to show that the limiting object of a stochastic process that represents the discretely updated weights fulfills \eqref{eq:evt}. For this purpose we define a stochastic process tracking the difference between both sides of \eqref{eq:evt}. More precisely, for} $\rho\in\MM$, $\eta\in\PP(\XX)$ and a stochastic process $(\theta_s)_s\sim\eta$, define a real-valued stochastic process by
\begin{equation}\label{eq:martingale_process}
\Phi_f^{(\rho,\eta)}(t, (\theta_s)_s) := f(\theta_t) - f(\theta_0) - \int_0^t \int_{\R^d} A(x, \theta_s; {\ms{Law}(\theta_s)}) \nabla f(\theta_s) \rho_s ({\rm d} x) {\rm d} s
\end{equation}
{for $0\le t\le T$, where $f\colon\R^d\to\R$ is a monomial of first order, i.e., given by
$$\R^d\ni x \mapsto x_i, \quad i\in\{1,\ldots,d\}$$}
and $\nabla$ represents the gradient. {We could now roughly follow \cite[Sections 4 and 5]{bufi} with a focus on \cite[Lemma 5.2]{bufi}, where a similar result is proven for a space of functions that are continuous with respect to the uniform topology and with continuous updates via a Brownian motion by invoking the connection between weak solutions to SDEs and local martingale problems to show that a limiting element of a sequence of (random) measures is indeed a solution to \eqref{eq:evt}. But in our case without a diffusion term, where the randomness in \eqref{eq:evt} comes only from the initial condition the arguments can be simplified. Roughly speaking, our situation is like a subcase of \cite[Lemma 3.3]{budconroy} in which there is a small common noise but no individual independent noises in particle dynamics. Indeed, in our situation the stochastic process $(\theta_s)_s\sim\eta$ is a solution to \eqref{eq:evt} if and only if for all monomials $f$ of first order, the process $\Phi_f^{(\rho,\eta)}$ is $\eta$-almost surely $0$ for all $0\le t\le T$.
Over the course of this section, we will show that if we plug the (random) limiting object from Lemma \ref{lem:weak_convergence_upper_bound} into \eqref{eq:martingale_process}, the respective stochastic process in \eqref{eq:martingale_process} is $\P$-almost surely $\bar\eta$-almost surely $0$ everywhere, which means it is a solution of \eqref{eq:evt}.}

Note that for any $\eta\in\PP(\XX)$, the process in \eqref{eq:martingale_process} is a {$\eta$-almost surely $0$ if for all $0\le t\le T$ and continuous bounded $\Psi\colon\XX\to\R$ it holds that
\begin{align}\label{eq:martingale_property}
    \begin{split}
    0 &= \E_{\eta}\big[\Psi(\cdot) \Phi_f^{(\rho,\eta)}(t, \cdot) \big] 
   := \int_{\XX} \Psi ((\theta_u)_u) \Phi_f^{(\rho,\eta)}(t, (\theta_u)_u) \eta({\rm d}(\theta_u)_u).
    \end{split}
\end{align}
This expression is now interpreted as a function with input $\rho$ and $\eta$. In \cite[Lemma 5.2]{bufi} a similar function as \eqref{eq:martingale_property} was essentially continuous with respect to $\rho$ and $\eta$ as part of the assumptions, which is not satisfied in our case.} Instead, we have to use truncations to induce boundedness and afterward approximate the limit with the truncated processes.
The following lemmas provide some help with this task. Recall from Section \ref{sec:mod} that $\CC$ denotes the set of continuous functions from $[0,T]$ to $\R^d$ endowed with the uniform topology. Further, recall that the sequence of measures $(\tilde\pi^n,\bar\eta^n)_n$ are constructed with common sequences of random distributions $((\pi_{k,n})_k)_n$ derived from Proposition \ref{thm:rep}. We additionally assume that the entropy bound from \eqref{eq:bound_entropy} or equivalently \eqref{eq:uniformly_bounded_entropy} is satisfied. In all of the following lemmas {$0\le t\le T$ can be any time point,} $f$ can be any monomial of first order from $\R^d$ to $\R$ and $\Psi\colon \XX \to \R$ any continuous bounded function.

\begin{lemma}[Mean convergence]\label{lem:mean_convergence}
Let $M>0$ and $(\rho^n,\eta^n)_n \subseteq \{\rho'\in\MM\colon \rho'\in\PP([0,T]\times[-M,M]^d)\}\times\PP(\{\theta\in\XX\colon \|\theta\|_\ff \le M\})$ converge weakly to some $(\rho,\eta) \in \{\rho'\in\MM\colon \rho'\in\PP([0,T]\times[-M,M]^d)\}\times\PP(\{\theta\in\XX\colon \|\theta\|_\ff \le M\})$. If $\eta(\CC) = 1$, it holds that
\begin{equation}\label{eq:convergence_martingale_process}
    {\E_{\eta^n}\big[\Psi(\cdot) \Phi_f^{(\rho^n,\eta^n)}(t, \cdot)\big] \overset{n\to\ff}{\longrightarrow} \E_{\eta}\big[\Psi(\cdot) \Phi_f^{(\rho,\eta)}(t, \cdot)\big].}
\end{equation}
\end{lemma}

\begin{lemma}[Integrable majorant]\label{lem:dominated_convergence}
Assume that $(\tilde\pi^n,\bar \eta^n)$ converges weakly to some $(\tilde\pi,\bar\eta)$ and that \ref{CON}, \ref{DCOMP} and \ref{WCOMP'} are satisfied. If almost surely $\bar\eta(\CC) = 1$, it holds that
\begin{align*}
    \E\bigg[&\int_\XX |f(\theta_t)| \bar\eta({\rm d}(\theta_u)_u) + \int_\XX |f(\theta_0)| \bar\eta({\rm d}(\theta_u)_u) \\
    &+ \int_{[0,t]\times \R^d} \int_{\XX} \|A(x, \theta_s; \bar\eta_s)\| \|\nabla f(\theta_s)\| \bar\eta({\rm d}(\theta_u)_u) \tilde\pi ({\rm d} (s,x))\bigg] < \ff.
\end{align*}
\end{lemma}

Next, for measures $\eta\in\PP(\XX)$ and $\rho\in\MM$ define push forwards $\eta|_M$ and $\rho|_M$ under the maps $\XX\ni\theta\mapsto (\theta \vee -M)\wedge M$ and $[0,T]\times\R^d\ni (t,x)\mapsto (t,(x\vee -M)\wedge M)$. The maximum or minimum of a vector and a real number is to be understood as a pointwise operation.

\begin{lemma}[Truncation approximation]\label{lem:weak_convergence_truncation}
Assume that $(\tilde\pi,\bar\eta)$ is the weak limit of $(\tilde\pi^n,\bar \eta^n)$, that $\bar\eta(\CC)=1$ and that \ref{CON}, \ref{DCOMP} and \ref{WCOMP'} are satisfied. Then,
\begin{align*}
    \E\Big[&\big|\E_{\bar\eta|_M}\big[\Psi(\cdot) \Phi_f^{(\tilde\pi|_M,\bar\eta|_M)}(t, \cdot)\big] - \E_{\bar\eta}\big[\Psi(\cdot) \Phi_f^{(\tilde\pi,\bar\eta)}(t, \cdot)\big]\big|\Big] \overset{M\to\ff}{\longrightarrow}0.
\end{align*}
\end{lemma}

\begin{lemma}[{Zero expecation}]\label{lem:almost_surely_zero}
Assume that \ref{CON} is satisfied. Then, for all $n\in\N$ and $M>0$ it holds that
\begin{equation}\label{eq:martingale_almost_surely_zero}
    {\E_{\bar\eta^n|_M}\big[ \Psi(\cdot) \Phi_f^{(\tilde\pi^n|_M,\bar\eta^n|_M)}(t, \cdot)\big] = 0.}
\end{equation}
\end{lemma}

\begin{proof}[Proof of Lemma \ref{lem:weak_convergence_upper_bound}]
In order to show that the weak limit $(\tilde\pi,\bar\eta)$ of $(\tilde\pi^n,\bar \eta^n)$ is indeed almost surely in $\PP_\ff^\nu$, we first deal with the initial distribution.
Since $\bar\eta_0^n$, the projection of $\bar\eta^n$ at time $0$ is deterministic and equal to $\nu_0^n$, we can deduce that $\bar\eta_0$ is equal to $\nu$ because of the assumption that $\nu_0^n$ converges weakly to $\nu$.

Next, we aim to apply Lemma \ref{lem:mean_convergence}. Recall that \eqref{eq:uniformly_bounded_entropy} implies that $\tilde\pi^n$ is almost surely absolute continuous with respect to $\pi_T$ for all $n\in\N$ because $\ff> \liminf_{n\to\ff} \E[R(\tilde\pi^n)] \ge \E[R(\tilde\pi)]$ by lower semicontinuity of $R$. Thus, $\tilde\pi\in\MM$.

To show that for almost all realizations of $\bar\eta$ it holds that $\bar\eta(\CC) = 1$, let $\Delta\co\XX\to[0,\ff)$, be the largest jump of a given function in $\XX$. Note that $\PP(\XX)$ is Polish as shown in the proof of Theorem \ref{thm:annealed}. The Skorokhod representation theorem \cite[Theorem 5.31]{kallenberg} yields the existence of random variables $(\hat\eta^n)_n$ and $\hat\eta$ in $\PP(\XX)$ defined on a common probability space $(\hat\Omega,\hat\AA,\hat\P)$ such that $\hat\eta^n$ converges $\hat\P$-almost surely to $\hat\eta$ and $\ms{Law}_\P(\bar\eta^n) = \ms{Law}_{\hat\P}(\hat\eta^n)$ as well as $\ms{Law}_\P(\bar\eta) = \ms{Law}_{\hat\P}(\hat\eta)$. By Lemma \ref{lem:max_bd}, for every $\e,\delta>0$ we can choose $n$ large enough such that
$$\P(\bar\eta^n(|\Delta|>\e)\le\delta) = 1.$$
Thus, we get
\begin{align*}
    \hat\P(\hat\eta(|\Delta|>\e)\le\delta) &\ge \hat\P\big(\liminf_{n\to\ff} \{\hat\eta^n(|\Delta|>\e)\le\delta\}\big) \ge \hat\P\big(\limsup_{n\to\ff} \{\hat\eta^n(|\Delta|>\e)\le \delta\}\big) \\
    &\ge \limsup_{n\to\ff} \hat\P(\hat\eta^n(|\Delta|>\e)\le \delta)= 1.
\end{align*}
Consequently, we get $\P(\bar\eta(|\Delta|>\e)\le\delta) = 1$ for all $\e,\delta>0$ and thus, $\P(\bar\eta(\CC) = 1) = 1$.

Now, we intend to prove that almost surely
\begin{equation}\label{eq:martingale_main_goal}
    {\E_{\bar\eta}\big[\Psi(\cdot) \Phi_f^{(\tilde\pi,\bar\eta)}(t, \cdot)\big] = 0.}
\end{equation}
For this, we consider the push forwards $\bar\eta^n|_M$, $\bar\eta|_M$, $\tilde\pi^n|_M$ and $\tilde\pi|_M$.
To these, we can apply Lemma \ref{lem:mean_convergence}, which states that if we restrict the domain accordingly, we have almost sure continuity with respect to the weak topology of the map
$$(\rho,\eta) \mapsto {\E_{\eta}\big[\Psi(\cdot) \Phi_f^{(\rho,\eta)}(t, \cdot) \big]}$$
for any monomial of first or second order $f\colon\R^d\to\R$, any ${0\le t\le T}$ and all continuous bounded $\Psi\colon \XX \to \R$.
Accordingly, the continuous mapping theorem yields
\begin{align*}
    &\lim_{n\to\ff} {\E_{\bar\eta^n|_M}\big[\Psi(\cdot) \Phi_f^{(\tilde\pi^n|_M,\bar\eta^n|_M)}(t, \cdot)\big] = \E_{\bar\eta|_M}\big[\Psi(\cdot) \Phi_f^{(\tilde\pi|_M,\bar\eta|_M)}(t, \cdot)\big]},
\end{align*}
where the convergence is to be understood as weak convergence.
From Lemma \ref{lem:weak_convergence_truncation} and Markov's inequality it follows that
\begin{align*}
    &{\E_{\bar\eta|_M}\big[\Psi(\cdot) \Phi_f^{(\tilde\pi|_M,\bar\eta|_M)}(t, \cdot)\big] \overset{M\to\ff}{\longrightarrow} \E_{\bar\eta}\big[\Psi(\cdot) \Phi_f^{(\tilde\pi,\bar\eta)}(t, \cdot)\big]},
\end{align*}
where this convergence is to be understood as convergence in probability. Thus, to show \eqref{eq:martingale_main_goal}, it suffices to show that for all $n\in\N$ and $M>0$ almost surely
$$
{\E_{\bar\eta^n|_M}\big[\Psi(\cdot) \Phi_f^{(\tilde\pi^n|_M,\bar\eta^n|_M)}(t, \cdot) \big] = 0}.
$$
But this is precisely the content of Lemma \ref{lem:almost_surely_zero}. Thus, we have shown \eqref{eq:martingale_main_goal} and deduce that the weak limit $(\tilde\pi,\bar\eta)$ of $(\tilde\pi^n,\bar \eta^n)$ is almost surely in $\PP_\ff^\nu$.
\end{proof}

We still need to prove Lemmas \ref{lem:mean_convergence}, \ref{lem:dominated_convergence}, \ref{lem:weak_convergence_truncation} and \ref{lem:almost_surely_zero}.
The proof of Lemma \ref{lem:mean_convergence} uses the induced boundedness, which enables us to invoke weak convergence arguments.

\begin{proof}[Proof of Lemma \ref{lem:mean_convergence}]
{When splitting up the expectation into its summable components, we get two parts that we can examine separately. First of all, we do not have to consider $f(\theta_0)$ separately and can instead only deal $f(\theta_t)$. The convergence of the part,}
\begin{equation}\label{eq:convergence_first_part_weak}
\int_{\XX} \Psi((\theta_u)_u) f(\theta_t) \eta^n({\rm d}(\theta_u)_u)
\end{equation}
follows with \cite[Theorem A.3.11]{dupuisellis}, which states that if the integrand is bounded and continuous at almost every point with respect to the weak limit of the measure, the integral converges. This can be applied since the integrand is bounded and continuous with respect to the uniform topology, $\eta$ has support in a subset of $\CC$ and because the Skorokhod topology coincides with the uniform topology, see \cite[page 124]{billingsley} on $\CC$. That $f$ is bounded here follows because $f$ is continuous and its argument is contained in the compact space $B_M(0)$.

Next, the convergence of the integration part of \eqref{eq:martingale_process} needs further justification. Thus, let us focus on the expectation with respect to $\eta^n$ of the integration part and rewrite it using Fubini's theorem to get
\begin{equation}\label{eq:convergence_second_part_weak}
\int_{[0,t]\times \R^d} \int_{\XX} \Psi((\theta_u)_u) A(x, \theta_r; \eta_r^n) \nabla f(\theta_r) \eta^n({\rm d}(\theta_u)_u) \rho^n ({\rm d} (r,x)).
\end{equation}
To show that this converges requires multiple steps. First, we can rely on \cite[pages 138, 139]{billingsley} to see that $\eta_r^n$ converges weakly to $\eta_r$ for all but countably many $r$. Fix an $r\in[0,t]$, where there is weak convergence. Now, an important observation is that we can write $A$ as
$$A(x,\theta_r;\eta_r^n) = g(x,\eta_r^n) h(x,\theta_r),$$
where $h \colon ((z,y),(c,w)) \mapsto (\s(z^\top w), c\s'(z^\top w)z)$ maps to $\R^d$ and both its arguments live in compact spaces. $g$ is by definition continuous and since its domain is a compact space, this makes $g$ uniformly continuous. This means that for any $\e>0$ we can choose $n$ large enough such that for all $x\in[-M,M]^d$ it holds that $|g(x,\eta_r^n)-g(x,\eta_r)| \le\e$. Thus, it holds for $n$ large enough that
\begin{align*}
&\Big\vert\int_{\R^d} \int_{\XX} \Psi((\theta_u)_u)g(x,\eta_r^n)-g(x,\eta_r) h(x,\theta_r) \nabla f(\theta_r) \eta^n({\rm d}(\theta_u)_u) \rho_r^n ({\rm d} x)\Big\vert \\
&\le\int_{\XX\times\R^d} \e |\Psi((\theta_u)_u) h(x,\theta_r) \nabla f(\theta_r)| (\eta^n\otimes\rho_r^n) ({\rm d}((\theta_u)_u, x)).
\end{align*}
With \cite[Theorem 2.8]{billingsley}, which states that product measures on separable spaces converge weakly if the marginal measures converge weakly, it follows that $\eta^n\otimes\rho_r^n$ converges weakly to $\eta\otimes\rho_r$. Thus, \cite[Theorem A.3.11]{dupuisellis} becomes applicable again and we can deduce that
\begin{align*}
&\int_{\XX\times\R^d} \e |\Psi((\theta_u)_u) h(x,\theta_r) \nabla f(\theta_r)| (\eta^n\otimes\rho_r^n) ({\rm d}((\theta_u)_u, x)) \\
&\overset{n\to\ff}{\longrightarrow} \int_{\XX\times\R^d} \e |\Psi((\theta_u)_u) h(x,\theta_r) \nabla f(\theta_r)| (\eta\otimes\rho_r) ({\rm d}((\theta_u)_u, x)) \overset{\e\to0}{\longrightarrow} 0,
\end{align*}
where the convergence to $0$ follows with the finiteness of the integral. The same argument yields that
\begin{align*}
&\lim_{n\to\ff} \int_{\R^d} \int_{\XX} \Psi((\theta_u)_u) g(x,\eta_r) h(x,\theta_r) \nabla f(\theta_r) \eta^n({\rm d}(\theta_u)_u) \rho^n ({\rm d} x) \\
&= \int_{\XX\times\R^d} \Psi((\theta_u)_u) g(x,\eta_r) h(x,\theta_r) \nabla f(\theta_r) (\eta\otimes\rho_r) ({\rm d}((\theta_u)_u, x)).
\end{align*}
This means that we get
\begin{align*}
    &\int_{\R^d} \int_\XX \Psi((\theta_u)_u) A(x, \theta_r; \eta_r^n) \nabla f(\theta_r) \eta^n({\rm d}(\theta_u)_u) \rho_r^n ({\rm d} x) \\
    &\overset{n\to\ff}{\longrightarrow} \int_{\R^d} \int_{\XX} \Psi((\theta_u)_u) A(x, \theta_r; \eta_r) \nabla f(\theta_r) \eta({\rm d}(\theta_u)_u) \rho_r ({\rm d} x)
\end{align*}
Finally, invoking dominated convergence again, we arrive at
\begin{align*}
&\int_{[0,t]} \int_{\R^d} \int_\XX \Psi((\theta_u)_u) A(x, \theta_r; \eta_r^n) \nabla f(\theta_r) \eta^n({\rm d}(\theta_u)_u) \rho_r^n ({\rm d} x) {\rm d} r \\
&= \int_{[0,t]\setminus\{t'\in[0,t]\colon\eta_{t'}^m\not\Rightarrow\eta_{t'}\}} \int_{\R^d} \int_\XX \Psi((\theta_u)_u) A(x, \theta_r; \eta_r^n) \nabla f(\theta_r) \eta^n({\rm d}(\theta_u)_u) \rho_r^n ({\rm d} x) {\rm d} r \\
&\overset{n\to\ff}{\longrightarrow} \int_{[0,t]\setminus\{t'\in[0,t]\colon\eta_{t'}^m\not\Rightarrow\eta_{t'}\}} \int_{\R^d} \int_\XX \Psi((\theta_u)_u) A(x, \theta_r; \eta_r) \nabla f(\theta_r) \eta({\rm d}(\theta_u)_u) \rho_r ({\rm d} x) {\rm d} r \\
&= \int_{[0,t]} \int_{\R^d} \int_\XX \Psi((\theta_u)_u) A(x, \theta_r; \eta_r) \nabla f(\theta_r) \eta({\rm d}(\theta_u)_u) \rho_r ({\rm d} x) {\rm d} r.
\end{align*}
This establishes the convergence of \eqref{eq:convergence_second_part_weak} and therefore, \eqref{eq:convergence_martingale_process}.
\end{proof}

The proof of Lemma \ref{lem:dominated_convergence} invokes weak convergence and monotone convergence to apply Lemma \ref{lem:th_gr}, which is used to show finiteness of the expectation.

\begin{proof}[Proof of Lemma \ref{lem:dominated_convergence}]
First, for any $\tilde d,M\in\N$ and $v\in\R^{\tilde d}$, we introduce the continuous map
$$\tau^M(v) := (|v_1|\wedge M,\dots,|v_{\tilde d}|\wedge M).$$
We start by considering the first summand in the statement of the lemma. First, by weak convergence, we get for all $M>0$
$$\E\bigg[\int_\XX \tau^M(f(\theta_t)) \bar\eta^n({\rm d}(\theta_u)_u)\bigg] \overset{n\to\ff}{\longrightarrow} \E\bigg[\int_\XX \tau^M(f(\theta_t)) \bar\eta({\rm d}(\theta_u)_u)\bigg].$$
Next, by monotone convergence, we get
$$\E\bigg[\int_\XX \tau^M(f(\theta_t)) \bar\eta({\rm d}(\theta_u)_u)\bigg] \overset{M\to\ff}{\longrightarrow} \E\bigg[\int_\XX |f(\theta_t)| \bar\eta({\rm d}(\theta_u)_u)\bigg].$$
Thus, there exist $C_1, C_2 >0$ and by Lemma \ref{lem:th_gr} (ii), $C_3 > 0$ such that
\begin{align*}
\E\bigg[\int_\XX |f(\theta_t)| \bar\eta({\rm d}(\theta_u)_u)\bigg] &\le \sup_M \E\bigg[\int_\XX \tau^M(f(\theta_t)) \bar\eta({\rm d}(\theta_u)_u)\bigg] + C_1 \le \sup_M \sup_n \E\bigg[\int_\XX \tau^M(f(\theta_t)) \bar\eta^n({\rm d}(\theta_u)_u)\bigg] + C_2 \\
&\le \sup_n \E\bigg[\int_\XX |f(\theta_t)| \bar\eta^n({\rm d}(\theta_u)_u)\bigg] + C_2 \le \sup_n \E\bigg[\int_\XX \|\theta_t\| \bar\eta^n({\rm d}(\theta_u)_u)\bigg] + C_2 \\
&\le C_3 \sup_n \E\big[(T^2 + 1) \big(\bar Y_n^{*,4} + \bar Z_n^{*,2}\big)\big] + C_3 <\ff.
\end{align*}
We used Cauchy-Schwarz to get the bound in the last line and the finiteness follows with Lemma \ref{lem:th_gr} (iii). {The second summand can be treated analogously.}

{To deal with the third summand}, we split it up into two parts that come from $g$, where we recall the definition of $A$ from \eqref{eq:grad}. {We point out that $\nabla f$ is of the form $(0,\dots,0,1,0,\dots,0)$ since $f$ is a monomial of first order. Thus, we do not need to truncate $\nabla f$ in the following.} We first get on the one hand that
\begin{align*}
    &\E\bigg[\int_{[0,t]\times \R^d} \int_{\XX} \bigg(\int_{\R^d} \tau^M(c') \bar\eta_s^n({\rm d}(c',w'))\bigg) \|(1, \tau^M(c_s) \tau^M(z))\| \cdot \|\nabla f(\theta_s)\| \bar\eta^n({\rm d}(c_u,w_u)_u) \tilde\pi^n ({\rm d} (s,x))\bigg] \\
    &\overset{n\to\ff}{\longrightarrow} \E\bigg[\int_{[0,t]\times \R^d} \int_{\XX} \bigg(\int_{\R^d} \tau^M(c') \bar\eta_s({\rm d}(c',w'))\bigg) \|(1, \tau^M(c_s) \tau^M(z)) \| \cdot \|\nabla f(\theta_s)\| \bar\eta({\rm d}(c_u,w_u)_u) \tilde\pi ({\rm d} (s,x))\bigg],
\end{align*}
when proceeding identical to the arguments after \eqref{eq:convergence_second_part_weak} and using that $\bar\eta(\CC)=1$ almost surely. Note that the pushforward of $\bar\eta^n$ given $\tau^M$ and the pushforward of $\tilde\pi^n$ given $\tau^M$ have compact support. On the other hand, monotone convergence again yields
\begin{align*}
    &\E\bigg[\int_{[0,t]\times \R^d} \int_{\XX} \bigg(\int_{\R^d} \tau^M(c') \bar\eta_s({\rm d}(c',w'))\bigg) \|(1, \tau^M(c_s) \tau^M(z))\| \cdot \|\nabla f(\theta_s)\| \bar\eta({\rm d}(c_u,w_u)_u) \tilde\pi ({\rm d} (s,x))\bigg] \\
    &\overset{M\to\ff}{\longrightarrow} \E\bigg[\int_{[0,t]\times \R^d} \int_{\XX} \bigg(\int_{\R^d} |c'| \bar\eta_s({\rm d}(c',w'))\bigg) \|(1, c_s z) \| \cdot \|\nabla f(\theta_s)\| \bar\eta({\rm d}(c_u,w_u)_u) \tilde\pi ({\rm d} (s,x))\bigg].
\end{align*}
We continue by computing for some $C_4, C_5, C_6>0$ that
\begin{align*}
    &\E\bigg[\int_{[0,t]\times \R^d} \int_{\XX} \bigg(\int_{\R^d} |c' \sigma(z^\top w')| \bar\eta_s({\rm d}(c',w'))\bigg) \|(\sigma(z^\top w_s), c_s \sigma' (z^\top w_s) z) \| \cdot \|\nabla f(\theta_s)\| \bar\eta({\rm d}(c_u,w_u)_u) \tilde\pi ({\rm d} (s,x))\bigg] \\
    &\le C_\s^2 \sup_M \E\bigg[\int_{[0,t]\times \R^d} \int_{\XX} \bigg(\int_{\R^d} \tau^M(c') \bar\eta_s({\rm d}(c',w'))\bigg) \|(1, \tau^M(c_s) \tau^M(z)) \| \\
    &\qquad\qquad\qquad\qquad\qquad\cdot \|\nabla f(\theta_s)\| \bar\eta({\rm d}(c_u,w_u)_u) \tilde\pi ({\rm d} (s,x))\bigg] + C_4\\
    &\le C_\s^2 \sup_M \sup_n \E\bigg[\int_{[0,t]\times \R^d} \int_{\XX} \bigg(\int_{\R^d} \tau^M(c') \bar\eta_s^n({\rm d}(c',w'))\bigg) \|(1, \tau^M(c_s) \tau^M(z))\| \\
    &\qquad\qquad\qquad\qquad\qquad\qquad\cdot \|\nabla f(\theta_s)\| \bar\eta^n({\rm d}(c_u,w_u)_u) \tilde\pi^n ({\rm d} (s,x))\bigg] + C_5 \\
    &\le C_\s^2 \sup_n \E\bigg[\int_{[0,t]\times \R^d} \int_{\XX} \bigg(\int_{\R^d} |c'| \bar\eta_s^n({\rm d}(c',w'))\bigg) \|(1, c_s z)\|  \cdot \|\nabla f(\theta_s)\| \bar\eta^n({\rm d}(c_u,w_u)_u) \tilde\pi^n ({\rm d} (s,x))\bigg] + C_5 \\
    &\le C_\s^2 \sup_n \E\bigg[(1+T\bar C (C_\nu + \bar Y_n^{*,1}))^2 \int_{[0,t]\times \R^d} \|(1, z)\| \tilde\pi^n ({\rm d} (s,x))\bigg] + C_5 \\
    &\le C_\s^2 \sup_n C_6 \E\bigg[\bar Y_n^{*,4} + \bigg(\int_{[0,t]\times \R^d} \|(1, z)\| \tilde\pi^n ({\rm d} (s,x))\bigg)^2\bigg] + C_6 =: (\star),
    \end{align*}
where we used \eqref{eq:bd_c} from the proof of Lemma \ref{lem:th_gr}, the Cauchy-Schwarz inequality, and {that $\|\nabla f\| \equiv 1$}. Now, we compute using the Cauchy-Schwarz inequality again that
\begin{align*}
    &\E\bigg[\bigg(\int_{[0,t]\times \R^d} \|(1,z)\| \tilde\pi^n ({\rm d} (s,x))\bigg)^2\bigg] \le T \E\bigg[\int_{[0,t]\times \R^d} \|(1,z)\|^2 \tilde\pi^n ({\rm d} (s,x))\bigg] \\
    &\le T \sum_{k=1}^{\lfloor nT\rfloor} \E\bigg[\int_{(k-1)/n}^{k/n} \int_{\R^d} \|(1,z)\|^2 \pi_{k,n} ({\rm d} x) {\rm d} s\bigg] = T \frac1n \sum_{k=1}^{\lfloor nT\rfloor} \E\bigg[\int_{\R^d} \|(1,z)\|^2 \pi_{k,n} ({\rm d} x) \bigg] \\
    &= T \frac1n \sum_{k=1}^{\lfloor nT\rfloor} \E\big[\E\big[ \|(1,\bar Z_{k,n})\|^2 \bigm\vert \FF_{k-1,n}\big]\big] = T \E[\bar Z_n^{*,2}].
\end{align*}
Now, Lemma \ref{lem:th_gr} (iii) implies that $(\star)<\ff$.

Further, we can conduct a similar calculation that yields
\begin{align*}
&\E\bigg[\int_{[0,t]\times \R^d} \int_{\XX} |y| \|(\sigma(z^\top w_s), c_s \sigma' (z^\top w_s) z) \|  \cdot \|\nabla f(\theta_s)\| \bar\eta({\rm d}(c_u,w_u)_u) \tilde\pi ({\rm d} (s,x))\bigg] < \ff.
\end{align*}
\end{proof}

The proof of Lemma \ref{lem:weak_convergence_truncation} invokes dominated convergence, which becomes applicable due to Lemma \ref{lem:dominated_convergence} to show that the truncation approximates the original.

\begin{proof}[Proof of Lemma \ref{lem:weak_convergence_truncation}]
To show the assertion note that for $x=(z,y)\in\R^d$, $(c,w)\in\R^d$ and $\eta\in\PP(\XX)$, we can write
$$A(x,(c,w),\eta) = y \alpha(z,(c,w)) - \int_{\R^d} \beta(z,(c',w')) \eta({\rm d}(c',w')) \alpha(z,(c,w)),$$
where $\alpha(z,(c,w)) := (\s(z^\top w), c\s'(z^\top w)z )$ and $\beta(z,(c,w)) := c \s(z^\top w)$. We also define $\tau^M(\,\cdot\,) := (\,\cdot\,\vee-M)\wedge M$, where the input can be a real number or a vector. In the latter case, the operations are to be understood pointwise. We use this to define
\begin{align*}
\alpha^M(z,(c,w)) &:= \big(\s(\tau^M(z)^\top \tau^M(w)), \tau^M(c)\s'(\tau^M(z)^\top \tau^M(w))\tau^M(z) \big), \\
\beta^M(z,(c,w)) &:= \tau^M(c) \s(\tau^M(z)^\top \tau^M(w))
\end{align*}
and for $(\theta_t)_t\in\XX$ set $\Psi^M((\theta_u)_u) := \Psi((\tau^M(\theta_u))_u)$. {We emphasize again that $\nabla f$ is a standard basis vector and to ease notation, in the next parts we will just write $\nabla f$ for the constant vector that this gradient represents.} Then, we can derive that
\begin{align*}
    &\E\bigg[\bigg| \int_{[{0},t]\times \R^d} \int_{\XX} \Psi^M((\theta_u)_u) \int_{\R^d} \beta^M(z,(c,w)) \bar\eta_r({\rm d}(c,w)) \alpha^M(z,\theta_r) \cdot\nabla f \bar\eta({\rm d}(\theta_u)_u) \tilde\pi ({\rm d} (r,(z,y))) \\
    &\quad- \int_{[{0},t]\times \R^d} \int_{\XX} \Psi((\theta_u)_u) \int_{\R^d} \beta(z,(c,w)) \bar\eta_r({\rm d}(c,w)) \alpha(z,\theta_r) \cdot\nabla f \bar\eta({\rm d}(\theta_u)_u) \tilde\pi ({\rm d} (r,(z,y))) \bigg| \bigg] \overset{M\to\ff}{\longrightarrow}0
\end{align*}
by dominated convergence where the dominating integrand is given by Lemma \ref{lem:dominated_convergence}. Note that since we only need to deal with monomials the gradient of $f$ is very simple. In the same manner, we get that
\begin{align*}
    &\E\bigg[ \bigg| \int_{[{0},t]\times \R^d} \int_{\XX} \Psi^M((\theta_u)_u) \tau^M(y) \alpha^M(z,\theta_r) \nabla f \bar\eta({\rm d}(\theta_u)_u) \tilde\pi ({\rm d} (r,(z,y))) \\
    &\quad- \int_{[{0},t]\times \R^d} \int_{\XX} \Psi((\theta_u)_u) y \alpha(z,\theta_r) \nabla f \bar\eta({\rm d}(\theta_u)_u) \tilde\pi ({\rm d} (r,(z,y))) \bigg| \bigg] \overset{M\to\ff}{\longrightarrow}0
\end{align*}
and
$$
    \E\bigg[ \bigg| \int_{\XX} \Psi^M((\theta_u)_u) f(\tau^M(\theta_t)) \bar\eta({\rm d}(\theta_u)_u) - \int_{\XX} \Psi((\theta_u)_u) f(\theta_t) \bar\eta({\rm d}(\theta_u)_u) \bigg|\bigg] \overset{M\to\ff}{\longrightarrow}0,
$$
{where we can also replace $t$ with $0$.}
\end{proof}

The proof of Lemma \ref{lem:almost_surely_zero} makes use of the dynamics of the update equation \eqref{eq:c2} to show that {the considered expectation is $0$}.

\begin{proof}[Proof of Lemma \ref{lem:almost_surely_zero}]
As in the proof of Lemma \ref{lem:weak_convergence_truncation} we define $\tau^M(\,\cdot\,) := (\,\cdot\vee-M)\wedge M$. With this, we set $A^M(x,(c,w),\eta) := \tau^M(y) \alpha^M(z,(c,w)) - \int_{\R^d} \beta^M(c',w') \eta({\rm d}(c',w')) \alpha^M(z,(c,w))$, where we reused the notation introduced in proof of Lemma \ref{lem:weak_convergence_truncation} and can rewrite \eqref{eq:martingale_almost_surely_zero} as
\begin{align}\label{eq:martingale_almost_surely_zero_rewrite}
    \begin{split}
    &\int_\XX \Psi^M((\theta_u)_u) \big(f(\tau^M(\theta_t)) - f(\tau^M(\theta_{{0}}))\big) \bar\eta^n({\rm d}(\theta_u)_u) \\
    &- \int_{[{0},t]\times \R^d} \int_{\XX} \Psi((\theta_u)_u) A^M(x,\theta_r; \bar\eta_r^n) \nabla f(\tau^M(\theta_r)) \bar\eta^n({\rm d}(\theta_u)_u) \tilde\pi^n ({\rm d} (r,x)).
    \end{split}
\end{align}
{Again we just write $\nabla f$ for the constant vector that this gradient represents.} The proof of Lemma \ref{lem:almost_surely_zero} is outlined by the following interim assertions.
\begin{enumerate}[leftmargin=*]
    \item[\textit{I.}] \emph{Rewriting:} We deal with the first summand in \eqref{eq:martingale_almost_surely_zero_rewrite} and show that almost surely
    $$\int_\XX \Psi^M((\theta_u)_u) \big(f(\tau^M(\theta_t)) - f(\tau^M(\theta_{{0}}))\big) \bar\eta^n({\rm d}(\theta_u)_u) = \Lambda_1^{f,n}({t}),$$
    where
    $$
        \Lambda_1^{f,n}({t}) := \frac1n \int_\XX \Psi^M((\theta_u)_u) \sum_{k = {1}}^{\lfloor nt\rfloor} A^M(\bar X_{k,n}, \theta_{(k-1)/n} ; \bar\eta_{(k-1)/n}^n) \nabla f \bar\eta^n({\rm d}(\theta_u)_u).
    $$
    \item[\textit{II.}] \emph{Approximation:} Next, we deal with the second summand of \eqref{eq:martingale_almost_surely_zero_rewrite} and show that almost surely
    \begin{align*}
        &\big\vert\int_{[{0},t]\times \R^d} \int_{\XX} \Psi^M((\theta_u)_u) A^M(x, \theta_r; \bar\eta_r^n) \cdot\nabla f \bar\eta^n({\rm d}(\theta_u)_u) \tilde\pi^n ({\rm d} (r,x)) - \Lambda_2^{f,n}({t})\big\vert \overset{n\to\ff}{\longrightarrow}0,
    \end{align*}
    where
    \begin{align*}
    &\Lambda_2^{f,n}({t}) := \frac1n \int_\XX \Psi^M((\theta_u)_u) \sum_{k = {1}}^{\lfloor nt\rfloor} \int_{\R^d} A^M(x, \theta_{(k-1)/n} ; \bar\eta_{(k-1)/n}^n) \cdot\nabla f \pi_{k,n}({\rm d} x) \bar\eta^n({\rm d}(\theta_u)_u).
    \end{align*}
    \item[\textit{III.}] \emph{Zero expectation step:} Finally, we show that the difference of $\Lambda_1^{f,n}$ and $\Lambda_2^{f,n}$ tends to $0$ in probability.
\end{enumerate}

Proof of \textit{I.}:

We can use the {simple representation of the gradient} of $f$. Note that for all appropriate $k,n$ it holds that $\theta_{k/n} = \theta_{(k-1)/n} + A(\bar X_{k,n}, \theta_{(k-1)/n} ; \bar\eta_{(k-1)/n}^n)/n$ if $(\theta_u)_u\sim\bar\eta^n$. More precisely, {by definition we} get almost surely
\begin{align*}
&\int_\XX \Psi^M((\theta_u)_u) \big(f(\tau^M(\theta_t)) - f(\tau^M(\theta_{{0}}))\big) \bar\eta^n({\rm d}(\theta_u)_u) \\
&= \int_\XX \Psi^M((\theta_u)_u) \sum_{k = {1}}^{\lfloor nt\rfloor} \big(f(\tau^M(\theta_{k/n})) - f(\tau^M(\theta_{(k-1)/n}))\big) \bar\eta^n({\rm d}(\theta_u)_u) = \Lambda_1^{f,n}({t})
\end{align*}

Proof of \textit{II.}:

The latter integral expression of \eqref{eq:martingale_almost_surely_zero_rewrite} can also be written as
\begin{align*}
&\int_\XX \Psi^M((\theta_u)_u) \int_{[{0},t]\times\R^d} A^M(x, \theta_r ; \bar\eta_r^n) \nabla f \tilde\pi^n({\rm d} r, {\rm d} x) \bar\eta^n({\rm d}(\theta_u)_u) \\
&= \underbrace{\int_\XX \Psi^M((\theta_u)_u) \sum_{k={1}}^{\lfloor nt\rfloor} \int_{(k-1)/n}^{k/n} \int_{\R^d} A^M(x, \theta_{(k-1)/n} ; \bar\eta_{(k-1)/n}^n) \cdot\nabla f \pi_{k,n}({\rm d} x) {\rm d} r \bar\eta^n({\rm d}(\theta_u)_u)}_{= \Lambda_2^{f,n}({t})} \\
&\quad+ \f1n \int_\XX \Psi^M((\theta_u)_u) \int_{t'(n)}^{t\wedge (\lfloor nT\rfloor/n)} \int_{\R^d} A^M(x, \theta_{\lfloor nt\rfloor/n} ; \bar\eta_{\lfloor nt\rfloor/n}^n)  \cdot \nabla f \pi_{\lfloor nt\rfloor+1,n}({\rm d} x) {\rm d} r \bar\eta^n({\rm d}(\theta_u)_u)
\end{align*}
for some $t'(n)\in[\lfloor nt\rfloor/n, t\wedge(\lfloor nT\rfloor/n)]$, where we point out that $\tilde\pi_{\lfloor nt\rfloor + 1,n}$ is not defined yet if $\lfloor nt\rfloor + 1>\lfloor nT\rfloor$. But in this case, the integration bounds of the respective integral collapse to a single point and thus, we can set $\tilde\pi_{\lfloor nt\rfloor + 1,n}:= \pi$. Due to the boundedness, we can deduce that
\begin{align}\label{eq:martingale2}
\begin{split}
&\Big|\int_{[{0},t]\times \R^d} \int_{\XX} \Psi^M((\theta_u)_u) A^M(x, \theta_r; \bar\eta_r^n) \nabla f \bar\eta^n({\rm d}(\theta_u)_u) \tilde\pi^n ({\rm d} (r,x)) - \Lambda_2^{f,n}({t})\Big| \\
&\le \f1n \Big|\int_\XX \Psi^M((\theta_u)_u) \int_{\R^d} A^M(x, \theta_{\lfloor nt\rfloor/n} ; \bar\eta_{\lfloor nt\rfloor/n}^n) \cdot \nabla f \pi_{\lfloor nt\rfloor+1,n}({\rm d} x) \bar\eta^n({\rm d}(\theta_u)_u)\Big| \overset{n\to0}{\longrightarrow} 0.
\end{split}
\end{align}

Proof of \textit{III.}:

In a first step, we can apply Cauchy-Schwarz inequality to get
\begin{align}\label{eq:martingale_diff}
\begin{split}
\E[(\Lambda_2^{f,n}({t}) - \Lambda_1^{f,n}({t}))^2] &\le \f1{n^2} \E\bigg[\Big(\int_\XX \Psi^M((\theta_u)_u) \bar\eta^n({\rm d}(\theta_u)_u)\Big)^2\bigg] \\
&\quad\cdot \E\bigg[\bigg(\sum_{k={1}}^{\lfloor nt\rfloor} \int_\XX \bigg(\int_{\R^d} A^M(x, \theta_{(k-1)/n} ; \bar\eta_{(k-1)/n}^n) \nabla f \pi_{k,n}({\rm d} x)\bigg) \\
&\qquad\qquad\qquad\qquad\quad- A^M(\bar X_{k,n}, \theta_{(k-1)/n} ; \bar\eta_{(k-1)/n}^n) \nabla f \bar\eta^n({\rm d}(\theta_u)_u)\bigg)^2\bigg]
\end{split}
\end{align}
Let $\bar\eta_{\le (k-1)/n}^n$ be the restriction of $\bar\eta^n$ to a random element in $\PP(D([0,(k-1)/n];\R^d))$ by not considering the behavior of functions in $\XX$ in the domain $((k-1)/n,T]$. By construction, $\bar\eta_{\le (k-1)/n}^n$ is $\FF_{k-1,n}$-measurable. After denoting $\XX_k:= D([0,(k-1)/n];\R^d)$, we compute
\begin{align*}
&\E\Big[\int_{\XX_k} A^M(\bar X_{k,n}, \theta_{(k-1)/n} ; \bar\eta_{(k-1)/n}^n) \nabla f \bar\eta_{\le (k-1)/n}^n({\rm d}(\theta_u)_u) \Bigm\vert \FF_{k-1,n}\Big] \\
&= \int_{\R^d} \int_{\XX_k} A^M(x, \theta_{(k-1)/n} ; \bar\eta_{(k-1)/n}^n) \nabla f \bar\eta_{\le (k-1)/n}^n({\rm d}(\theta_u)_u) \tilde\pi_{k,n}({\rm d} x) \\
&= \int_{\XX_k} \int_{\R^d} A^M(x, \theta_{(k-1)/n} ; \bar\eta_{(k-1)/n}^n) \nabla f \tilde\pi_{k,n}({\rm d} x) \bar\eta_{\le (k-1)/n}^n({\rm d}(\theta_u)_u).
\end{align*}
Since also $\tilde\pi_{k,n}$ is $\FF_{k-1,n}$-measurable, we have
\begin{align*}
&\E\Big[\int_\XX \int_{\R^d} A^M(x, \theta_{(k-1)/n} ; \bar\eta_{(k-1)/n}^n) \nabla f \tilde\pi_{k,n}({\rm d} x) - A^M(\bar X_{k,n}, \theta_{(k-1)/n} ; \bar\eta_{(k-1)/n}^n) \nabla f \bar\eta^n({\rm d}(\theta_u)_u) \Big] \\
&= \E\Big[\int_{\XX_k} \int_{\R^d} A^M(x, \theta_{(k-1)/n} ; \bar\eta_{(k-1)/n}^n) \nabla f \tilde\pi_{k,n}({\rm d} x) \bar\eta_{\le (k-1)/n}^n({\rm d}(\theta_u)_u) \\
&\qquad\quad\ - \E\Big[\int_{\XX_k} A^M(\bar X_{k,n}, \theta_{(k-1)/n} ; \bar\eta_{(k-1)/n}^n) \cdot \nabla f \bar\eta_{\le (k-1)/n}^n({\rm d}(\theta_u)_u) \Bigm\vert \FF_{k-1,n}\Big]\Big] = 0.
\end{align*}
Consequently, the sequence
\begin{align*}
    &\bigg(\int_{\R^d} (A^M(x, \theta_{(k-1)/n} ; \bar\eta_{(k-1)/n}^n) - A^M(\bar X_{k,n}, \theta_{(k-1)/n} ; \bar\eta_{(k-1)/n}^n)) \nabla f \tilde\pi_{k,n}({\rm d} x)\bigg)_k
\end{align*}
is a martingale difference sequence and we deduce that when expanding the square in \eqref{eq:martingale_diff} the cross-terms disappear. Thus,
\begin{align*}
\E[(\Lambda_2^{f,n}({t}) - \Lambda_1^{f,n}({t}))^2] &\le \f1{n^2} \E\bigg[\Big(\int_\XX \Psi^M((\theta_u)_u) \bar\eta^n({\rm d}(\theta_u)_u)\Big)^2\bigg] \\
&\quad\cdot \sum_{k={1}}^{\lfloor nt\rfloor} \E\bigg[\bigg( \int_\XX \bigg(\int_{\R^d} A^M(x, \theta_{(k-1)/n} ; \bar\eta_{(k-1)/n}^n) \nabla f \tilde\pi_{k,n}({\rm d} x)\bigg) \\
&\qquad\qquad\qquad\ \ - A^M(\bar X_{k,n}, \theta_{(k-1)/n} ; \bar\eta_{(k-1)/n}^n) \nabla f \bar\eta^n({\rm d}(\theta_u)_u)\bigg)^2\bigg] \overset{n\to\ff}{\longrightarrow} 0,
\end{align*}
where the convergence follows since all integrands are bounded and thus, we multiply $1/n^2$ with something of order $n$.
\end{proof}

\subsection{Proofs -- lower bound}

{Now, we prove} Lemmas \ref{lem:weak_convergence_lower_bound} and \ref{lem:entropy_bound}. The former involves some integral calculation combined with weak convergence, while for the latter we use the Donsker-Varadhan variational formula \cite[Lemma 2.4 (a)]{dupuis} to derive the entropy bound.

\begin{proof}[Proof of Lemma \ref{lem:weak_convergence_lower_bound}]
Let $f\co\R^d\to\R$ and $g\co[0,T]\to \R$ be bounded and continuous.
If we can show convergence of $\int_0^T\int_{\R^d} f(x)g(t) \rho_t^n({\rm d} x){\rm d} t$ to $\int_0^T\int_{\R^d} f(x)g(t) \rho_t({\rm d} x){\rm d} t$, the weak convergence follows with \cite[Theorem A.3.14]{dupuisellis}, which states it suffices to consider the product of two bounded continuous functions when determining weak convergence on the product space of two Polish spaces.

	In order to show this, for $k\le nT$, we abbreviate $I_{k,n} := [(k-1)/n,k/n]$, and with Fubini {and $T_n :=\lfloor nT\rfloor/n $,} we can derive
\begin{align*}
`\int_0^{T_n}\int_{\R^d} f(x)g(s) \rho_s^n({\rm d} x){\rm d} s &= n \sum_{k\le nT} \int_0^{T_n} \int_{I_{k,n}} \int_{\R^d} f(x) g(s) \one\{s\in I_{k,n}\} \rho_t({\rm d} x) {\rm d} t {\rm d} s \\
&= n \sum_{k\le nT} \int_{I_{k,n}} \int_{\R^d} \int_{I_{k,n}} f(x) g(s) {\rm d} s \rho_t({\rm d} x) {\rm d} t \\
&= n\sum_{k\le nT} \int_0^{T_n} \int_{\R^d} \int_{I_{k,n}} f(x) g(s) \one\{t\in I_{k,n}\} {\rm d} s \rho_t({\rm d} x) {\rm d} t \\
&= \int_0^{T_n} \int_{\R^d} n \int_{I_{\lfloor nt\rfloor,n}} f(x) g(s) {\rm d} s \rho_t({\rm d} x) {\rm d} t.
\end{align*}
Next, the previous calculation yields
\begin{align*}
&\Big\vert \int_0^T\int_{\R^d} f(x)g(s) \rho_s^n({\rm d} x){\rm d} s - \int_0^T\int_{\R^d} f(x)g(t) \rho_t({\rm d} x){\rm d} t \Big\vert \\
&= \Big\vert \int_0^{\lfloor nT\rfloor/n} \int_{\R^d} n \int_{I_{\lfloor nt\rfloor,n}} f(x) g(s) {\rm d} s \rho_t({\rm d} x) {\rm d} t - \int_0^{\lfloor nT\rfloor/n} \int_{\R^d} f(x)g(t) \rho_t({\rm d} x){\rm d} t \Big\vert \\
&\quad+ \Big\vert\int_{\lfloor nT\rfloor/n}^T \int_{\R^d} f(x)g(t) \rho_t({\rm d} x){\rm d} t\Big\vert \\
&\le \int_0^{\lfloor nT\rfloor/n} \int_{\R^d} n \int_{I_{\lfloor nt\rfloor,n}} \vert f(x) g(s) - f(x) g(t)\vert {\rm d} s \rho_t({\rm d} x) {\rm d} t + (T-\lfloor nT\rfloor /n) \sup_{t\in[0,T], x\in\R^d} |f(x)g(t)|.
\end{align*}
Since $g$ is continuous on a compact interval, it is uniformly continuous and for an arbitrary $\e>0$ we can choose $n$ large enough such that if $|s-t|\le 1/n$ it also holds that $|g(s)-g(t)|\le\e$. Therefore, for $n$ large enough we arrive at
\begin{align*}
    &\Big\vert \int_0^T\int_{\R^d} f(x)g(s) \rho_s^n({\rm d} x){\rm d} s - \int_0^T\int_{\R^d} f(x)g(t) \rho_t({\rm d} x){\rm d} t \Big\vert \\
    &\le T \e \sup_{x\in\R^d} |f(x)| + (T-\lfloor nT\rfloor /n) \sup_{t\in[0,T], x\in\R^d} |f(x)g(t)|.
\end{align*}
	{Since $f$ and $g$ are bounded, the assertion follows, when we let $n\to\ff$ and $\e\to0$.}
\end{proof}

\begin{proof}[Proof of Lemma \ref{lem:entropy_bound}]
Let $\CC_b(\R^d)$ denote the space of bounded continuous maps from $\R^d$ to $\R$. Then, we compute with Fubini
\begin{align*}
R(\rho^n) &= \int_0^T H(\rho_s^n \mid \pi) {\rm d} s = \int_0^T \sup_{g\in\CC_b(\R^d)} \bigg( \int_{\R^d} g(x) \rho_s^n({\rm d} x) - \log\int_{\R^d} e^{g(x)} \pi({\rm d} x) \bigg) {\rm d} s \\
&= \int_0^T \sup_{g\in\CC_b(\R^d)} \bigg( \int_{\R^d} g(x) \sum_{k\le nT} \one\{s\in[(k-1)/n,k/n]\} \int_{(k-1)/n}^{k/n} n \rho_t({\rm d} x){\rm d} t - \log\int_{\R^d} e^{g(x)} \pi({\rm d} x) \bigg) {\rm d} s \\
&= \int_0^T \sup_{g\in\CC_b(\R^d)} \bigg( n \sum_{k\le nT} \one\{s\in[(k-1)/n,k/n]\} \cdot\int_{(k-1)/n}^{k/n} \bigg(\int_{\R^d} g(x) \rho_t({\rm d} x) - \log\int_{\R^d} e^{g(x)} \pi({\rm d} x) \bigg) {\rm d} t\bigg) {\rm d} s.
\end{align*}
Therefore,
\begin{align*}
R(\rho^n)&\le n \sum_{k\le nT} \int_0^T \one\{s\in[(k-1)/n,k/n]\} \cdot\int_{(k-1)/n}^{k/n} \sup_{g\in\CC_b(\R^d)} \bigg(\int_{\R^d} g(x) \rho_t({\rm d} x) - \log\int_{\R^d} e^{g(x)} \pi({\rm d} x) \bigg) {\rm d} t {\rm d} s \\
&= \int_0^T \sup_{g\in\CC_b(\R^d)} \bigg(\int_{\R^d} g(x) \rho_t({\rm d} x) - \log\int_{\R^d} e^{g(x)} \pi({\rm d} x) \bigg) {\rm d} t = R(\rho),
\end{align*}
where to get the first and last equality we applied \cite[Corollary 2.7]{dupuis}, the chain rule for the relative entropy, to deal with the product of measures. The second and second-to-last equality follows with the Donsker-Varadhan variational formula, which can be found for instance in \cite[Lemma 2.4 (a)]{dupuis}.
\end{proof}

\subsection{Proofs -- good rate function}
To prove Lemma \ref{lem:compact_level_sets} we can follow the proof of \cite[Proposition 8.1]{kill} which also proves goodness of a rate function. However, the rate function that we derived is completely different from the one considered in \cite{kill}, which is also used as a tightness function in their proof. This is not possible with our more complicated rate function and thus, we need to proceed differently to establish tightness. For this reason, we insert an additional lemma that is helpful here and in the proof of Corollary \ref{cor:weakLLN}.

\begin{lemma}[Bounded entropy implies compactly supported trajectories]\label{lem:compact_support_bounded_entropy}
Assume \ref{CON}, \ref{DCOMP} and \ref{WCOMP'} and let $M>0$. Then, there exists a compact set $\bar K$ such that all ${\eta}\in\PP(\XX)$ with $I_\nu({\eta})\le M$ are supported on $\bar K$.
\end{lemma}

The proof of Lemma \ref{lem:compact_level_sets} requires us to {revisit the characterization of solutions of \eqref{eq:evt} with the process in \eqref{eq:martingale_process} that is $0$ mentioned in Section \ref{sec:laplace_upper_bound}. Recall the notation of $\Phi_f^{(\rho,\eta)}$ from \eqref{eq:martingale_process} for $\rho\in\MM$, $\eta\in\PP(\XX)$ and monomials of first order $f$. This is $\eta$-almost surely $0$ if and only if for all $0\le t\le T$ and continuous bounded $\Psi\colon\XX\to\R$ equation \eqref{eq:martingale_property} holds.
Our aim is to}
apply Lemma \ref{lem:mean_convergence} again for which we need slightly modified versions of Lemma \ref{lem:dominated_convergence} and \ref{lem:weak_convergence_truncation} to perform the truncation argument. This is necessary because Lemma \ref{lem:dominated_convergence} and \ref{lem:weak_convergence_truncation} deal with limiting objects originating from the measure construction via SGD, while now we are interested in measures that are already contained in a space of admissible measures specified in Section \ref{sec:mod} and are not necessarily the limits of measures constructed via SGD. For this, we remind the reader of more of the notation introduced before Lemma \ref{lem:mean_convergence} and \ref{lem:weak_convergence_truncation} and state the following lemma.

\begin{lemma}[Integrable majorant and truncation approximation]\label{lem:dominated_convergence_and_truncation_approx}
Assume \ref{CON} and \ref{DCOMP}. Let $\eta\in\PP(\XX)$ be supported on a compact set $\bar K$ and let $\rho\in\MM$ such that $R(\rho)<\ff$. Then, for all $t\in[0,T]$ it holds that
\begin{align*}
    &\int_\XX |f(\theta_t)| \eta({\rm d}(\theta_u)_u) + \int_\XX |f(\theta_0)| \eta({\rm d}(\theta_u)_u) + \int_{[0,t]\times \R^d} \int_{\XX} \|A(x, \theta_s; \eta_s)\| \|\nabla f(\theta_s)\| \eta({\rm d}(\theta_u)_u) \rho ({\rm d} (s,x)) < \ff.
\end{align*}
Further, it holds that
$$
    {|\E_{\eta}\big[\Psi(\cdot) \Phi_f^{(\rho|_N,\eta)}(t, \cdot) \big] - \E_{\eta}\big[\Psi(\cdot) \Phi_f^{(\rho,\eta)}(t, \cdot)\big]\big|} \overset{N\to\ff}{\longrightarrow}0.
$$
\end{lemma}

\begin{proof}[Proof of Lemma \ref{lem:compact_level_sets}]
For an arbitrary $M>0$, let $(\nu'_n,{\eta^n})_n\subseteq\{(\nu',{\eta})\in\PP(\{(c,w)\in\R^d: \|(c,w)\|\le C_\nu\})\times\PP(\XX)\colon I_{\nu'}({\eta})\le M\}$. Since
$$M\ge \sup_n I_{\nu'_n}({\eta^n}) = \sup_n \inf_{{(\rho,\tilde\eta)} \in \PP_\ff^{\nu'_n}\co \tilde\eta=\eta^n} \frac1{T} R({\rho}),$$
there exists a sequence $(\rho^n)_n\subseteq\MM$ such that $(\rho^n,{\eta^n})\in\PP_\ff^{\nu'_n}$ and $\frac1{T} R(\rho^n) \le M+\frac1n$ for each $n\in\N$. Since $R$ has compact sublevel sets, this means that $(\rho^n)_n$ is contained in a compact set, which makes it tight by Prokhorov's theorem. Further, by Lemma \ref{lem:compact_support_bounded_entropy}, $({\eta^n})_n$ have a uniform compact support, also implying tightness.
Besides that, $\PP(\{(c,w)\in\R^d: \|(c,w)\|\le C_\nu\})$ is compact, which makes $(\nu'_n)_n$ tight. Let $(\nu',\rho,{\eta})$ be a limit point of a subsequence of $(\nu'_n,\rho^n,{\eta^n})_n$. To avoid overcomplicating the notation, we keep identifying this subsequence by $n$. First, by lower semicontinuity of the relative entropy, it holds that
$\frac1{T} R(\rho) \le \liminf_{n\to\ff} \frac1{T} R(\rho^n) \le M.$

	Next, we assert that $(\rho,{\eta})$ corresponds to a weak solution to \eqref{eq:evt} with respect to $\nu'$ as initial distribution, i.e., $(\theta_t)_t\sim{\eta}$ is a weak solution to \eqref{eq:evt}. In particular, we assert that $(\rho,{\eta})\in\PP_\ff^{\nu'}$. This implies that $I_{\nu'}({\eta}) \le \frac1{T} R(\rho)\le M$, meaning that $(\nu',{\eta})$ is in the set of the statement. Then, the Eberlein-Šmulian theorem \cite[Theorem 2.8.6]{megginson} yields compactness and finalizes the proof.

To show the assertion note that since ${\eta^n}(\CC) = 1$ for all $n$ by the definition of $\PP_\ff^{\nu_n'}$, also ${\eta}(\CC)=1$ has to be satisfied as follows for instance from \cite[Theorem 13.4]{billingsley}, which states that the weak limit of measures on a Skorokhod space with decreasing jump sizes is supported on the space of continuous functions with respect to the uniform topology. Further, by Lemma \ref{lem:compact_support_bounded_entropy}, the $({\eta^n})_n$ have uniformly bounded support. This enables us to apply Lemma \ref{lem:mean_convergence} and Lemma \ref{lem:dominated_convergence_and_truncation_approx}. The latter lets us approximate the data trajectory measure with $\rho|_N$, the push forward of $\rho$ under a truncation for some $N\in\N$, which can then be handled by Lemma \ref{lem:mean_convergence}. Together, similar as in the proof of Lemma \ref{lem:weak_convergence_upper_bound}, this yields that
$\lim_{n\to\ff} {\E_{{\eta^n}}\big[\Psi(\cdot) \Phi_f^{(\rho^n,\eta^n)}(t, \cdot)\big] = \E_{\eta}\big[\Psi(\cdot) \Phi_f^{(\rho,\eta)}(t, \cdot)\big]}.$
Since $(\rho^n,{\eta^n})\in\PP_\ff^{\nu'_n}$, i.e., $(\theta_t)_t\sim{\eta^n}$ is a weak solution of \eqref{eq:evt} with respect to the initial weight distribution $\nu'_n$, {the stochastic process $\Phi_f^{(\rho^n,\eta^n)}$ is $\eta^n$-almost surely $0$ as outlined in the beginning of Section \ref{sec:laplace_upper_bound}, therefore,} the left side is equal to zero and thus, also the right side is zero. Finally, it has to hold that the projection ${\eta_0^n}$ of ${\eta^n}$ at time $0$ is equal to $\nu'_n$. Otherwise, $I_{\nu'_n}({\eta^n}) =\ff$, since there could not exist an {$\eta$ that together with some data trajectory measure is part of $\PP_\ff^{\nu'_n}$} such that ${\eta}={\eta^n}$. Because $\nu'_n$ converges weakly to  $\nu'$ this implies that ${\eta_0^n}$ converges weakly to $\nu'$, which then has to be equal to ${\eta_0}$. This proves that $(\rho,{\eta})\in\PP_\ff^{\nu'}$.
\end{proof}

\begin{proof}[Proof of Lemma \ref{lem:compact_support_bounded_entropy}]
	We showed the {Laplace} lower bound in Section \ref{sec:proof_overview_lower_bound}. From this, the lower bound of the LDP follows without the property that $I_\nu$ has to be a good rate function, see \cite[Proof of Theorem 1.8]{dupuis}. The good rate function property is only required to derive the upper large deviations bound from the upper bound of the Laplace principle. Therefore, we can use the lower large deviations bound to get for any open (with respect to the weak topology) set $O\subseteq \PP(\XX)$ that
\begin{equation}\label{eq:large_deviations_lower_bound}
\liminf_{n\to\ff}\frac1{n'} \log\P(\eta^n\in O) \ge -\inf_{\eta\in O} I_\nu(\eta).
\end{equation}
Next, let $\delta>0$ and $t_1,t_2\in[0,T]$ with $\delta< t_2-t_1\le 2\delta$ as well as $L>0$. We intend to show that for the particular open sets $\{\eta\in\PP(\XX) \colon \eta(\|\omega_{t_2} - \omega_{t_1}\| > \e) > 0\}^\circ$, where $^\circ$ denotes the interior (with respect to the weak topology) and $\{\eta\in\PP(\XX) \colon \eta(\|\omega\|_\ff > L) > 0\}$ the left side of \eqref{eq:large_deviations_lower_bound} becomes arbitrarily small with decreasing $\delta$ and increasing $L$, respectively. That means we assert that
\begin{equation}\label{eq:large_deviations_lower_bound_1}
    \limsup_{\delta\to 0} \sup_{\substack{t_1,t_2\in [0,T]\colon \\ \delta<t_2-t_1\le 2\delta}} \liminf_{n\to\ff} \frac1{n'} \log\P\big(\eta^n \in\{\eta\in\PP(\XX) \colon \eta(\|\omega_{t_2} - \omega_{t_1}\| > \e) > 0\}\big) = -\ff
\end{equation}
and
\begin{equation}\label{eq:large_deviations_lower_bound_2}
    \limsup_{L\to\ff} \liminf_{n\to\ff} \frac1{n'} \log\P\big(\eta^n \in\{\eta\in\PP(\XX) \colon \eta(\|\omega\|_\ff > L) > 0\}\big) = -\ff.
\end{equation}
Note that if $I_\nu({\eta}) < M$, then ${\eta}(\CC)=1$ by the definition of $\PP_\nu^\ff$. Since $\{{\tilde\eta}\in\PP(\CC) \colon {\tilde\eta}(\|\omega_{t_2} - \omega_{t_1}\| > \e) > 0\} \subseteq \{{\tilde\eta}\in\PP(\XX) \colon {\tilde\eta}(\|\omega_{t_2} - \omega_{t_1}\| > \e) > 0\}$ is an open set, we can deduce from \eqref{eq:large_deviations_lower_bound} and \eqref{eq:large_deviations_lower_bound_1} that if $I_\nu({\eta}) < M$ it is also satisfied that ${\eta}(\|\omega_{t_2} - \omega_{t_1}\| > \e) = 0$ for all $\e>0$ and small enough $\delta>0$.
Thus,
$$\sup_{t_1,t_2\in [0,T]\colon \delta<t_2-t_1\le 2\delta} {\eta}(\|\omega_{t_2} - \omega_{t_1}\| \ge \e) = 0.$$
Similarly, we get that for ${\eta}\in\PP(\XX)$ with $I_\nu({\eta}) < M$ it holds that ${\eta}(\|\omega\|_\ff > L) = 0$ for large enough $L>0$ by \eqref{eq:large_deviations_lower_bound} and \eqref{eq:large_deviations_lower_bound_2}. This means that we can find a compact set $\bar K$ such that all ${\eta}$ with $I_\nu({\eta})\le M$ are supported on a compact set by the characterization of compactness in \cite[Theorem 13.2]{billingsley}.

	It remains to prove the assertions \eqref{eq:large_deviations_lower_bound_1} and \eqref{eq:large_deviations_lower_bound_2}. We start with the first one and since \ref{CON} and \ref{WCOMP'} were assumed, we can {reuse} the bound in \eqref{eq:update_bound} from the proof of Lemma \ref{lem:th_gr}. If we recall $(\hat\theta_t^n)_t$ from \eqref{eq:th_def} and set $({\bar\theta_t^n})_t := (\hat\theta_{\lfloor tn\rfloor/n}^n((X_{k,n})_{k\le \lfloor tn\rfloor}))_t$, now plugging in the unaltered data points instead of the tilted ones and set $Y_n^* := \frac1n\sum_{k=1}^{\lfloor nT\rfloor} |Y_{k,n}|$, we get that for every $\omega_0\in B_{C_\nu}(0)$
\begin{align*}
    &\frac1n \sum_{k=\lfloor nt_1\rfloor + 1}^{\lfloor nt_2\rfloor} \|A(X_{k,n}, {\bar\theta_{(k-1)/n}^n}; \eta_{(k-1)/n}^n)\| \\
    &\le \frac1n \sum_{k=\lfloor nt_1\rfloor + 1}^{\lfloor nt_2\rfloor}\big(|Y_{k,n}| + C_\s^2 C_\nu e^{C_\s^2 T} + C_\s^3 e^{C_\s^2 T} Y_n^* \big) \bar C \big(C_\nu + Y_n^*\big) \|(1, Z_{k,n})\| \\
    &\le \frac{C_1}{n} \sum_{k=\lfloor nt_1\rfloor + 1}^{\lfloor nt_2\rfloor} \Big(\big(|Y_{k,n}| \|(1, Z_{k,n})\| + \|(1, Z_{k,n})\| + \|(1, Z_{k,n})\| Y_n^* + |Y_{k,n}| \|(1, Z_{k,n})\| Y_n^* + \|(1, Z_{k,n})\| (Y_n^*)^2 \big)\Big)
\end{align*}
for some $C_1>0$.
Next, we introduce the notation $\mu_1 := \log \E[e^{|Y_{1,n}|}]$ and $\mu_2(\e) := \log \E[e^{\e^{-4} \|(1,Z_{1,n})\|}]$ and use exponential Markov's inequality to get 
\begin{align*}
    \P\bigg(\frac1n \sum_{k=\lfloor nt_1\rfloor + 1}^{\lfloor nt_2\rfloor}\|(1, Z_{k,n})\| (Y_n^*)^2 >\e\bigg)
    &\le \P\bigg(\frac1n \sum_{k=\lfloor nt_1\rfloor + 1}^{\lfloor nt_2\rfloor}\|(1, Z_{k,n})\| (Y_n^*)^2 > \e, Y_n^*\ge \e^{-1}\bigg) \\
    &\quad+ \P\bigg(\frac1n \sum_{k=\lfloor nt_1\rfloor + 1}^{\lfloor nt_2\rfloor}\|(1, Z_{k,n})\| (Y_n^*)^2 > \e, Y_n^* < \e^{-1}\bigg) \\
    &\le \P\big(Y_n^* \ge \e^{-1}\big) + \P\bigg(\frac1n \sum_{k=\lfloor nt_1\rfloor + 1}^{\lfloor nt_2\rfloor}\|(1, Z_{k,n})\| > \e^3\bigg) \\
    &\le \P\big(Y_n^* \ge \e^{-1}\big) + \P\bigg(\frac{\e^{-4}}{n} \sum_{k=\lfloor nt_1\rfloor + 1}^{\lfloor nt_2\rfloor}\|(1, Z_{k,n})\| > \e^{-1}\bigg) \\
    &\le \E[e^{|Y_{1,n}|}]^{nT} e^{-\e^{-1} n} + \E[e^{\e^{-4} \|(1,Z_{1,n})\|}]^{2\delta n + 1} e^{-\e^{-1} n} \\
    &\le e^{-(\e^{-1}-T\mu_1)n} + e^{-(-\mu_2(\e)(2\delta + 1/n) + \e^{-1})n}.
\end{align*}
Note that $\mu_1$ and $\mu_2(\e)$ are finite by \ref{DCOMP} and that the same calculations yield similar bounds for the other summands. Now, let $0<\tilde\e<\e$. Using an $\tilde\e/5$ argument together with the union bound we get for some $C_2>0$
\begin{align*}
    &\limsup_{\delta\to 0} \sup_{\substack{t_1,t_2\in [0,T]\colon \\ \delta<t_2-t_1\le 2\delta}} \limsup_{n\to\ff} \frac1{n'} \log \sup_{\omega_0\in B_{C_\nu}(0)} \P\big( \|{\bar\theta_{t_2}^n} - {\bar\theta_{t_1}^n}\| > \e\big) \\
    &\le \limsup_{\delta\to 0} \sup_{\substack{t_1,t_2\in [0,T]\colon \\ \delta<t_2-t_1\le 2\delta}} \limsup_{n\to\ff} \frac1{n'} \log \sup_{\omega_0\in B_{C_\nu}(0)} \P\big( \|{\bar\theta_{t_2}^n} - {\bar\theta_{t_1}^n}\| > \tilde\e\big) \\
    &\le \limsup_{\delta\to 0} \sup_{\substack{t_1,t_2\in [0,T]\colon \\ \delta<t_2-t_1\le 2\delta}} \limsup_{n\to\ff}\frac1{n'}  \cdot\log \bigg(\sup_{\omega_0\in B_{C_\nu}(0)} \P\bigg( \frac1n \sum_{k=\lfloor nt_1\rfloor + 1}^{\lfloor nt_2\rfloor} \|A(X_{k,n}, {\bar\theta_{(k-1)/n}^n}; \eta_{(k-1)/n}^n)\| > \tilde\e\bigg)\bigg) \\
    &\le -C_2 \tilde\e^{-1} \overset{\tilde\e\to 0}{\longrightarrow} -\ff.
\end{align*}
Finally, this implies that
\begin{align*}
    &\limsup_{\delta\to 0} \sup_{\substack{t_1,t_2\in [0,T]\colon \\ \delta<t_2-t_1\le 2\delta}} \liminf_{n\to\ff} \frac1{n'} \log\P\big(\eta^n \in\{\eta\in\PP(\XX) \colon \eta(\|\omega_{t_2} - \omega_{t_1}\| > \e) > 0\}\big) \\
    &\le \limsup_{\delta\to 0} \sup_{\substack{t_1,t_2\in [0,T]\colon \\ \delta<t_2-t_1\le 2\delta}} \limsup_{n\to\ff} \frac1{n'} \log \sup_{\omega_0\in B_{C_\nu}(0)} \P\big(\|{\bar\theta_{t_2}^n} - {\bar\theta_{t_1}^n}\| > \e\big) = -\ff,
\end{align*}
which proves \eqref{eq:large_deviations_lower_bound_1}. The proof of \eqref{eq:large_deviations_lower_bound_2} is analogous. More, precisely, we use the bound from Lemma \ref{lem:th_gr} (ii) to get
\begin{align*}
    \sup_{\omega_0\in B_{C_\nu}(0)} \P\big( \|({\bar\theta_t^n})_t\|_\ff > L\big) \le \P\big(C_\nu + C_\ms{SGD}(T^2 + 1) \big(Y_n^{*,4} + Z_n^{*,2}\big) > L\big).
\end{align*}
Next, let $\mu_3 := \log \E[e^{|Y_{1,n}|^4}]$ and use exponential Markov's inequality to get for one of the summands
\begin{align*}
\P\big( Y_n^{*,4} > n(L - C_\nu)/(2(T^2 + 1) C_\ms{SGD}) \big) &\le \E[e^{|Y_{1,n}|^4}]^{nT} e^{-n(L - C_\nu)/(2 (T^2 + 1) C_\ms{SGD})} \\
&\le e^{-n(-T\mu_3 + (L - C_\nu)/(2 (T^2 + 1) C_\ms{SGD}))}.
\end{align*}
Repeating the same steps from above we can conclude \eqref{eq:large_deviations_lower_bound_2}.
\end{proof}

\begin{proof}[Proof of Lemma \ref{lem:dominated_convergence_and_truncation_approx}]
Due to \ref{CON} and since $\eta$ is supported on $\bar K$ it is sufficient to prove
$$\int_{[0,t]\times \R^d} |y| + \|(1,z)\| \rho ({\rm d} (s,x)) < \ff.$$
This can be shown by proceeding as in the proof of Lemma \ref{lem:th_gr} (iii). Using the inequality $a b \le e^a + (b\log b - b + 1)$ for $a,b\ge 0$ from \cite[Equation (2.9)]{dupuis} we compute
\begin{align*}
	\int_{[0,t]\times\R^d} |y| \rho({\rm d} (z,y)) &\le \int_{[0,t]\times\R^d} |y| \frac{{\rm d}\rho}{{\rm d}\pi_T}(s,x) \pi_T({\rm d} (s,x)) \le \int_{[0,t]\times\R^d}e^{|y|} \pi_T({\rm d} (s,x)) + R(\rho) < \ff
\end{align*}
where the finiteness follows with \ref{DCOMP}. The second term containing $\|(1,z)\|$ can be handled similarly. The rest of the proof a simplified version of the proof of Lemma \ref{lem:weak_convergence_truncation}, which we omit here.
\end{proof}

\section{Proofs of Theorem \ref{thm:unq} and Corollary \ref{cor:weakLLN}}\label{sec:unq}

This section is dedicated to proving Theorem \ref{thm:unq}, i.e., that having the additional assumption of a compactly supported data distribution is sufficient for \ref{UNQ} and proving Corollary \ref{cor:weakLLN}, the weak LLN. We start with the latter.

\begin{proof}[Proof of Corollary \ref{cor:weakLLN}]
First, note that given $\pi_T$ as data trajectory distribution and $\nu$ as initial weight distribution, \eqref{eq:evt} has a weakly unique solution ${\eta}$ by \ref{UNQ}.
Now, let $\AA\subseteq\PP(\XX)$ be closed with ${\eta}\not\in\AA$ and assume $\inf_{{\eta'}\in \AA} J({\eta'}) = 0$. Then, by the definition of $J$, there exist sequences $(\nu^n)_n$ and $(\rho^n,{\tilde\eta^n})\in\PP_\ff^{\nu^n}$ for each $n\in\N$ with ${\tilde\eta^n}\in\AA$ such that
$$H(\nu_0^n\mid\nu) + R(\rho^n) \le \frac1n.$$
In particular, this means that $H(\nu_0^n\mid\nu)\le 1$ as well as $R(\rho^n) \le 1$. Thus, by Lemma \ref{lem:tr} and since $\nu$ has compact support, both sequences are tight. Further, by Lemma \ref{lem:compact_support_bounded_entropy}, all ${\tilde\eta^n}$ share a uniform compact support and are therefore also tight. Denote by $\nu^*$, $\rho^*$ and ${\tilde\eta}$ subsequential limits of $(\nu^n)_n$, $(\rho^n)_n$ and ${\tilde\eta^n}$, respectively such that $(\nu^n,\rho^n,{\tilde\eta^n})$ converges weakly to $(\nu^*,\rho^*,{\tilde\eta})$ along the same subsequence. For notational convenience, we keep identifying this subsequence by $n$. Note that due to lower semicontinuity of the relative entropy, we have that
$$H(\nu^*\mid\nu) \le \liminf_{n\to\ff} H(\nu_0^n\mid\nu) \le \liminf_{n\to\ff} 1/n = 0,$$
showing that $\nu^* = \nu$. The same arguments alongside the same calculation show that $\rho^* = \pi_T$.
Now, using Lemma \ref{lem:mean_convergence} and Lemma \ref{lem:dominated_convergence_and_truncation_approx} we can argue analogous to the proof of Lemma \ref{lem:compact_level_sets} to deduce that $(\pi_T,{\tilde\eta})\in\PP_\ff^{\nu}$. Additionally, \ref{UNQ} implies that ${\tilde\eta} = {\eta}$. Thus, ${\eta}\in\AA$, which is a contradiction. Hence, $\inf_{{\eta'}\in\AA} J({\eta'}) > 0$.

In particular, for any bounded and continuous function $f\colon\XX\to\R$ and $\e>0$ we can use Theorem \ref{thm:annealed} to deduce that
\begin{align*}
    \limsup_{n\to\ff} \P(|\theta^n(f) - {\eta}(f)| \ge \e) &= \limsup_{n\to\ff} \P(\theta^n \in \underbrace{\{{\eta'}\in\PP(\XX)\colon |{\eta'}(f) - {\eta}(f)| \ge \e\}}_{=:\AA}) \le \limsup_{n\to\ff} e^{-n \inf_{{\eta'}\in\AA} J({\eta'})} = 0.
\end{align*}
\end{proof}

To prove Theorem \ref{thm:unq}, we proceed as in \cite[Section 4]{siri}. However, our setting is more involved since we have to deal with tilted measures of $\pi_T$ instead.
We insert two additional helpful lemmas for the proof of Theorem \ref{thm:unq}.

\begin{lemma}[Bounded trajectories]\label{lem:bounded_trajectory}
If \ref{WCOMP} holds and the data distribution $\pi$ has compact support, then a weak solution of \eqref{eq:evt} connected to an element $\PP_\ff^\nu$ is bounded, i.e., there is $C_\ms{traj} > 0$ such that the weak solution $(\theta_t)_t$ satisfies $\|(\theta_t)_t\|_\ff \le C_\ms{traj}$.
\end{lemma}

\begin{proof}
Let $(\rho,\eta)\in\PP_\ff^\nu$. Since the data distribution $\pi$ was assumed to have compact support, the definition of $\PP_\ff^\nu$ implies that the support of $\rho$ is a subset of $\supp(\pi)$.
Now, by the definition of the gradient in \eqref{eq:grad} and \eqref{eq:evt} together with \ref{WCOMP} yield the boundedness of the weak solution connected to $\eta$.
\end{proof}

Before we state the second lemma, for any $\rho\in\MM$ with $R(\rho)<\ff$, we define a map $\zeta^\rho\co \PP(\XX) \to \PP(\XX)$ similar to \cite[Section 4]{siri} such that for $\eta \in\PP(\XX)$ it holds that $\zeta^\rho(\eta) = \ms{Law}((\theta_t)_t)$, where $(\theta_t)_t$ is a stochastic process given by
\begin{align*}
\theta_t = \theta_0 + \int_0^t \int_{\R^d} A(x,\theta_s; \eta_s) \rho_s ({\rm d} x) {\rm d} s, \quad \theta_0 \sim \nu.
\end{align*}
If we split $(\theta_t)_t$ into its weights $(c_t)_t$ associated to the output layer and the weights $(w_t)_t$ associated by the hidden layer, this means
\begin{align}\label{eq_initial_value_problem}
\begin{split}
&c_t = c_0 + \int_0^t \int_{\R^d}  (y - F(z, \eta_s)) \s(z^\top w_s) \rho_s ({\rm d} (z,y)) {\rm d} s \\
&w_t = w_0 + \int_0^t \int_{\R^d} (y - F(z, \eta_s)) c_s \s'(z^\top w_s)z \rho_s ({\rm d} (z,y)) {\rm d} s \\
&(c_0,w_0) \sim \nu.
\end{split}
\end{align}
In accordance with Lemma \ref{lem:bounded_trajectory} and the continuity of weak solutions that are connected to an element in $\PP_\ff^\nu$, we can restrict ourselves to measures in $\PP(\bar\CC(T))$, where $\bar\CC(T) := \{(\theta_t)_t\in\CC(T) \colon \|(\theta_t)_t\|_\ff \le C_\ms{traj}\}$, where for any $0\le T_0\le T$, the set $\CC(T_0)$ are functions that are continuous with respect to the uniform topology from $[0,T_0]$ to $\R^d$. Further, for any $T_0\in[0,T]$, let $\zeta_{T_0}^\rho$ be defined as $\zeta^\rho$ but restricted to the domain $\PP(\bar\CC(T_0))$.

Additionally, for any $T_0\in[0, T]$, we define the Wasserstein metric on $\PP(\bar\CC(T_0))$ as in \cite{siri}. This means that, for $\eta,\eta'\in\PP(\bar\CC(T_0))$ we set
$$\DD_{T_0}(\eta,\eta') := \inf_{\mu\in\PP_{\eta,\eta'}(\bar\CC(T_0)\times\bar\CC(T_0))} \bigg(\int_{\bar\CC(T_0)\times\bar\CC(T_0)} \sup_{s\le T_0} \|x_s-y_s\|^2 \mu {\rm d}(x,y)\bigg)^{1/2},$$
where $\PP_{\eta,\eta'}(\bar\CC(T_0)\times\bar\CC(T_0))$ are probability measures on the product space $\bar\CC(T_0)\times\bar\CC(T_0)$ with marginal distributions corresponding to $\eta$ and $\eta'$. Now, we proceed as in \cite[Section 4]{siri}.

\begin{lemma}[Contraction property]\label{lem:contraction_property}
Assume \ref{CON}, \ref{WCOMP} and that $\pi$ has compact support. Let $\rho\in\MM$ with $R(\rho)<\ff$. Then, there exists $C_\ms{contr}>0$ such that for all $T_0\in[0,T]$ and $\eta^1,\eta^2\in \PP(\bar\CC(T_0))$ it holds that
$$\DD_{T_0}(\zeta_{T_0}^\rho(\eta^1),\zeta_{T_0}^\rho(\eta^2)) \le T_0 C_\ms{contr} \DD_{T_0}(\eta^1,\eta^2).$$
\end{lemma}

\begin{proof}[Proof of Theorem \ref{thm:unq}]
Due to \ref{WCOMP} and the compactly supported distribution $\pi$, Lemma \ref{lem:bounded_trajectory} becomes applicable, which with the definition of $\PP_\ff^\nu$ implies that we can indeed restrict ourselves to weak solutions that are stochastic processes valued in $\bar\CC(T)$.
Note that by definition, if such a stochastic process $(\theta_t)_t$ in $\bar\CC(T)$ is a weak solution of \eqref{eq:evt}, its law $\ms{Law}((\theta_t)_t)$ has to be a fixed point of $\zeta^\rho$ and vice versa if $\eta\in\PP(\bar\CC(T))$ is a fixed point of $\zeta^\rho$ then $(\theta_t)_t\sim\eta$ is a weak solution of \eqref{eq:evt}. Thus, it suffices to show that $\zeta^\rho$ has a unique fixed point. To this end, we invoke the Banach fixed-point theorem and restrict ourselves to functions defined on $[0,T_0]$ for some $T_0\in[0,T]$.

By the Stone-Weierstrass theorem, $\CC(T_0)$ equipped with the supremum norm is separable and if we further restrict the function space to functions bounded by $C_\ms{traj}$, then it stays separable and is also still complete and it induces boundedness of the Wasserstein metric. Then, we can invoke the main theorem of \cite{bolley2008} to see that $\DD_{T_0}$ turns $\PP(\bar\CC(T_0))$ into a Banach space. This is a condition for the Banach fixed-point theorem. Further, since we also assumed \ref{CON}, Lemma \ref{lem:contraction_property} becomes applicable. Now, by Lemma \ref{lem:contraction_property}, when choosing $T_0$ sufficiently small we can turn $\zeta_{T_0}^\rho$ into a contraction.
Finally, this lets us follow the rest of the arguments from \cite[Proof of Lemma 4.4]{siri}. More precisely, partitioning $[0, T]$ into intervals of length at most $T_0$ and repeating this argument proves uniqueness and existence of a fixed point of $\zeta^\rho$ in $\PP(\bar\CC(T))$ and therefore weak uniqueness of \ref{eq:evt}.
\end{proof}

Lemma \ref{lem:contraction_property} can be shown by replicating the proof of \cite[Lemma 4.3]{siri}.

\begin{proof}[Proof of Lemma \ref{lem:contraction_property}]
First, for $\eta^1,\eta^2\in \PP(\bar\CC(T_0))$, let $(c_t^1,w_t^1)_t$ and $(c_t^2,w_t^2)_t$ be distributed according to $\zeta_{T_0}(\eta^1)$ and $\zeta_{T_0}(\eta^2)$, respectively. This means that $(c_t^1,w_t^1)_t$ and $(c_t^2,w_t^2)_t$ satisfy \eqref{eq_initial_value_problem}, which lets us compute
\begin{align}\label{eq:contraction_property}
\begin{split}
c_{T_0}^1-c_{T_0}^2 &= \int_0^{T_0} \int_{\R^d}  (y - F(z, \eta_s^1)) \s(z^\top w_s^1) \rho_s ({\rm d} (z,y)) {\rm d} s - \int_0^{T_0} \int_{\R^d}  (y - F(z, \eta_s^2)) \s(z^\top w_s^2) \rho_s ({\rm d} (z,y)) {\rm d} s + c_0^1 - c_0^2 \\
&= \int_0^{T_0} \int_{\R^d} y (\s(z^\top w_s^1) - \s(z^\top w_s^2)) \rho_s ({\rm d} (z,y)) {\rm d} s + \int_0^{T_0} \int_{\R^d}  F(z, \eta_s^2) (\s(z^\top w_s^1) - \s(z^\top w_s^2)) \rho_s ({\rm d} (z,y)) {\rm d} s \\
&\quad+ \int_0^{T_0} \int_{\R^d} (F(z, \eta_s^2) - F(z,\eta_s^1)) \s(z^\top w_s^1) \rho_s ({\rm d} (z,y)) {\rm d} s + c_0^1 - c_0^2.
\end{split}
\end{align}
We examine the first of the terms. We assumed $\R(\rho)<\ff$. Thus, $\supp(\rho)\subseteq\supp(\pi)\subseteq B_{C_\pi}(0)$ for some $C_\pi\ge 1$. Using this and \ref{CON} we obtain
\begin{align*}
&\Big\vert\int_0^{T_0} \int_{\R^d}  y (\s(z^\top w_s^1) - \s(z^\top w_s^2)) \rho_s ({\rm d} (z,y)) {\rm d} s\Big\vert \le \int_0^{T_0} \int_{\R^d} |y| L_\s |z^\top (w_s^1 - w_s^2)| \rho_s ({\rm d} (z,y)) {\rm d} s \\
&\le \int_0^{T_0} \int_{\R^d} |y| L_\s \|z\| \|w_s^1 - w_s^2\| \rho_s ({\rm d} (z,y)) {\rm d} s \le L_\s C_\pi^2 \int_0^{T_0} \|w_s^1 - w_s^2\| {\rm d} s,
\end{align*}
where we applied the Cauchy-Schwarz inequality to get the second inequality. We get a similar bound for the second term by invoking that $\supp(\eta^2)\subseteq\PP(\bar\CC(T_0))$ and thus, $|F(z,\eta_s^2)| \le C_\s C_\ms{traj}$, which lets us calculate
\begin{align*}
\Big\vert\int_0^{T_0} \int_{\R^d} F(z, \eta_s^2) (\s(z^\top w_s^1) - \s(z^\top w_s^2)) \rho_s ({\rm d} (z,y)) {\rm d} s\Big\vert
&\le \int_0^{T_0} \int_{\R^d} |F(z, \eta_s^2)| L_\s \|z\| \|w_s^1 - w_s^2\| \rho_s ({\rm d} (z,y)) {\rm d} s \\
&\le C_\ms{traj} C_\s L_\s C_\pi \int_0^{T_0} \|w_s^1 - w_s^2\| {\rm d} s.
\end{align*}
To deal with the last term of \eqref{eq:contraction_property}, for any $\mu\in\PP_{\eta^1,\eta^2}(\bar\CC(T_0),\bar\CC(T_0))$ note that
\begin{align*}
|F(z, \eta_s^2) - F(z,\eta_s^1)|
&= \Big\vert\int_{\R^d\times\R^d} c \s(z^\top w) - c' \s(z^\top w') \mu_s({\rm d}(c',w',c,w))\Big\vert \\
&\le \int_{\R^d\times\R^d} |(c-c') \s(z^\top w) + c'( \s(z^\top w) -\s(z^\top w'))| \mu_s({\rm d}(c',w',c,w)) \\
&\le \int_{\R^d\times\R^d} C_\s |c-c'| +  C_\ms{traj} L_\s C_\pi \|w -w'\| \mu_s({\rm d}(c',w',c,w)).
\end{align*}
This leads to
\begin{align*}
&\Big\vert\int_0^{T_0} \int_{\R^d} (F(z, \eta_s^2) - F(z,\eta_s^1)) \s(z^\top w_s^1) \rho_s ({\rm d} (z,y)) {\rm d} s\Big\vert \\
&\le \int_0^{T_0} \int_{\R^d\times\R^d} C_\s^2 |c-c'| +  C_\s C_\ms{traj} L_\s C_\pi \|w -w'\| \mu_s({\rm d}(c',w',c,w)) \rho_s ({\rm d} (z,y)) {\rm d} s \\
&\le C_\s^2 C_\ms{traj} L_\s C_\pi \int_0^{T_0} \bigg(\int_{\R^d\times\R^d} (|c-c'| + \|w -w'\|)^2 \mu_s({\rm d}(c',w',c,w))\bigg)^{1/2} {\rm d} s \\
&\le C_\s^2 C_\ms{traj} L_\s C_\pi \int_0^{T_0} \bigg(\int_{\R^d\times\R^d} 2(|c-c'|^2 + \|w -w'\|^2) \mu_s({\rm d}(c',w',c,w))\bigg)^{1/2} {\rm d} s \\
\end{align*}
where we used Cauchy-Schwarz inequality again and utilized the inequality $(a+b)^2 \le 2 (a^2 + b^2)$ for $a,b\in\R$. Further, after transitioning from the distribution of the projection $\mu_s$ to $\mu$ we arrive at
\begin{align*}
&\int_0^{T_0} \bigg(\int_{\R^d\times\R^d} (|c-c'|^2 + \|w -w'\|^2) \mu_s({\rm d}(c',w',c,w))\bigg)^{1/2} {\rm d} s \\
&= \int_0^{T_0} \bigg(\int_{\bar\CC(T_0)\times\bar\CC(T_0)} (|c_s-c'_s|^2 + \|w_s -w'_s\|^2) \mu({\rm d}((c'_t,w'_t)_t,(c_t,w_t)_t))\bigg)^{1/2} {\rm d} s \\
&\le T_0 \bigg(\int_{\bar\CC(T_0)\times\bar\CC(T_0)} \sup_{u\le T_0} (|c_u-c'_u|^2 + \|w_u -w'_u\|^2) \mu({\rm d}((c'_t,w'_t)_t,(c_t,w_t)_t))\bigg)^{1/2}.
\end{align*}
Since, $\mu\in\PP_{\eta^1,\eta^2}(\bar\CC(T_0),\bar\CC(T_0))$ was chosen arbitrarily, it follows that
$$\Big\vert\int_0^{T_0} \int_{\R^d} (F(z, \eta_s^2) - F(z,\eta_s^1)) \s(z^\top w_s^1) \rho_s ({\rm d} (z,y)) {\rm d} s\Big\vert \le T_0 2 C_\s^2 C_\ms{traj} L_\s C_\pi \DD_{T_0}(\eta^1,\eta^2).$$
Together we get
$$|c_{T_0}^1-c_{T_0}^2| \le C_1 \int_0^{T_0} \|w_s^1-w_s^2\| {\rm d} s+ T_0 C_1 \DD_{T_0}(\eta^1,\eta^2) + |c_0^1 - c_0^2|$$
for $C_1 := 2 C_\s^2 C_\ms{traj} L_\s C_\pi^2$ and with an analogous computation that is adapted from \cite[Appendix B]{siri} we can examine $w^1_{T_0}-w^2_{T_0}$ and obtain
\begin{align*}
w_{T_0}^1-w_{T_0}^2 &= \int_0^{T_0} \int_{\R^d} (y - F(z, \eta_s^1)) c_s^1 \s'(z^\top w_s^1) z \rho_s ({\rm d} x) {\rm d} s \\
&\quad- \int_0^{T_0} \int_{\R^d} (y - F(z, \eta_s^2)) c_s^2 \s'(z^\top w_s^2) z \rho_s ({\rm d} x) {\rm d} s + w_0^1 - w_0^2 \\
&= \int_0^{T_0} \int_{\R^d} y z (c_s^1 \s'(z^\top w_s^1) - c_s^2 \s'(z^\top w_s^2)) \rho_s ({\rm d} x) {\rm d} s \\
&\quad+ \int_0^{T_0} \int_{\R^d} F(z, \eta_s^2) z (c_s^2 \s'(z^\top w_s^2) - c_s^1 \s'(z^\top w_s^1)) \rho_s ({\rm d} x) {\rm d} s \\
&\quad+ \int_0^{T_0} \int_{\R^d} (F(z, \eta_s^2) - F(z,\eta_s^1)) z c_s^1 \s'(z^\top w_s^1) \rho_s ({\rm d} x) {\rm d} s + w_0^1 - w_0^2. \\
&= \int_0^{T_0} \int_{\R^d} y z (c_s^1 - c_s^2) \s'(z^\top w_s^2) \rho_s ({\rm d} x) {\rm d} s + \int_0^{T_0} \int_{\R^d} y z c_s^1 (\s'(z^\top w_s^1) - \s'(z^\top w_s^2)) \rho_s ({\rm d} x) {\rm d} s \\
&\quad+ \int_0^{T_0} \int_{\R^d} F(z, \eta_s^2) z (c_s^2 - c_s^1) \s'(z^\top w_s^2) \rho_s ({\rm d} x) {\rm d} s \\
&\quad+ \int_0^{T_0} \int_{\R^d} F(z, \eta_s^2) z c_s^1 (\s'(z^\top w_s^2) - \s'(z^\top w_s^1)) \rho_s ({\rm d} x) {\rm d} s \\
&\quad+ \int_0^{T_0} \int_{\R^d} (F(z, \eta_s^2) - F(z,\eta_s^1)) z c_s^1 \s'(z^\top w_s^1) \rho_s ({\rm d} x) {\rm d} s + w_0^1 - w_0^2.
\end{align*}
In this case, we get the following bounds for the summands. For the first term, we have
\begin{align*}
    \Big\vert\int_0^{T_0} \int_{\R^d} y z (c_s^1 - c_s^2) \s'(z^\top w_s^2) \rho_s ({\rm d} x) {\rm d} s \Big\vert &\le \int_0^{T_0} \int_{\R^d} |y| \|z\| |c_s^1 - c_s^2| C_\s \rho_s ({\rm d} x) {\rm d} s \le C_\s C_\pi^2 \int_0^{T_0} |c_s^1 - c_s^2| {\rm d} s,
\end{align*}
for the second one
\begin{align*}
&\Big\vert\int_0^{T_0} \int_{\R^d} y z c_s^1 (\s'(z^\top w_s^1) - \s'(z^\top w_s^2)) \rho_s ({\rm d} x) {\rm d} s\Big\vert \\
&\le \int_0^{T_0} \int_{\R^d} |y| \|z\|^2 L_\s C_\ms{traj} \|w_s^1 - w_s^2\| \rho_s ({\rm d} x) {\rm d} s \le C_\pi^3 L_\s C_\ms{traj} \int_0^{T_0} \|w_s^1 - w_s^2\| {\rm d} s,
\end{align*}
for the third one
\begin{align*}
&\Big\vert\int_0^{T_0} \int_{\R^d} F(z, \eta_s^2) z (c_s^2-c_s^1) \s'(z^\top w_s^2) \rho_s ({\rm d} x) {\rm d} s\Big\vert \\
&\le \int_0^{T_0} \int_{\R^d} C_\ms{traj} C_\s^2 \|z\| |c_s^1-c_s^2| \rho_s ({\rm d} x) {\rm d} s \le C_\pi  C_\ms{traj} C_\s^2 \int_0^{T_0} |c_s^1-c_s^2| {\rm d} s,
\end{align*}
for the fourth one
\begin{align*}
&\Big\vert \int_0^{T_0} \int_{\R^d} F(z, \eta_s^2) z c_s^1 (\s'(z^\top w_s^2) - \s'(z^\top w_s^1)) \rho_s ({\rm d} x) {\rm d} s\Big\vert \\
&\le \int_0^{T_0} \int_{\R^d} C_\ms{traj}^2 C_\s \|z\|^2 L_\s \|w_s^1 - w_s^2\| \rho_s ({\rm d} x) {\rm d} s \le C_\ms{traj}^2 C_\s C_\pi^2 L_\s \int_0^{T_0} \|w_s^1 - w_s^2\| {\rm d} s
\end{align*}
and for the fifth one, the previous calculation does not change despite having the additional factors $z$ and $c_s^1$ in the integrand, which give an additional factor $C_\ms{traj} C_\pi$ in the bound. Finally, combining these bounds and setting $C_2 := 2 C_\ms{traj}^2 C_\pi^3 C_\s^2 L_\s$, we get
$$\|w_{T_0}^1-w_{T_0}^2\| \le C_2 \int_0^{T_0} |c_s^1-c_s^2| + \|w_s^1-w_s^2\| {\rm d} s + {T_0} C_2 \DD_{T_0}(\eta^1,\eta^2) + \|w_0^1 - w_0^2\|.$$
Jointly expressed this means
\begin{align*}
    |c_{T_0}^1-c_{T_0}^2| + \|w_{T_0}^1-w_{T_0}^2\| &\le C \int_0^{T_0} |c_s^1-c_s^2| + \|w_s^1-w_s^2\| {\rm d} s + {T_0} C \DD_{T_0}(\eta^1,\eta^2) {\rm d} s + |c_0^1 - c_0^2| + \|w_0^1 - w_0^2\|
\end{align*}
for $C := C_1+C_2$, which implies
\begin{align*}
    \sup_{s\le {T_0}} \big(|c_s^1-c_s^2| + \|w_s^1-w_s^2\|\big) &\le C \int_0^{T_0} \sup_{u\le s} \big(|c_s^1-c_s^2| + \|w_s^1-w_s^2\|\big) {\rm d} s + {T_0} C \DD_{T_0}(\eta^1,\eta^2) + |c_0^1 - c_0^2| + \|w_0^1 - w_0^2\|.
\end{align*}
Here, with Gronwall's lemma, we obtain
\begin{align*}
    &\sup_{s\le T_0} \big(|c_s^1-c_s^2|^2 + \|w_s^1-w_s^2\|^2\big) \le \Big(\sup_{s\le T_0} \big(|c_s^1-c_s^2| + \|w_s^1-w_s^2\|\big)\Big)^2 \\
    &\le \big(T_0 C \DD_{T_0}(\eta^1,\eta^2) + |c_0^1 - c_0^2| + \|w_0^1 - w_0^2\|\big)^2 e^{2C T_0}
\end{align*}
from which follows that for any $\mu\in \PP_{\zeta_{T_0}^\rho(\eta^1),\zeta_{T_0}^\rho(\eta^2)}(\bar\CC(T_0),\bar\CC(T_0))$ we can calculate with the abbreviation $D^{(0)}(\mu) :=|c_0^1 - c_0^2| + \|w_0^1 - w_0^2\|$
\begin{align*}
&\bigg(\int_{\bar\CC(T_0)\times\bar\CC(T_0)} \sup_{s\le T_0} \|(c_s^1,w_s^1)-(c_s^2,w_s^2)\|^2 \mu({\rm d}((c_t^1,w_t^1)_t,(c_t^2,w_t^2)_t))\bigg)^{1/2} \\
&\le e^{C T_0} \bigg(\big(T_0 C \DD_{T_0}(\eta^1,\eta^2) \big)^2 + 2 T_0 C \DD_{T_0}(\eta^1,\eta^2) \int_{\bar\CC(T_0)\times\bar\CC(T_0)} D^{(0)}(\mu) \mu({\rm d}((c_t^1,w_t^1)_t,(c_t^2,w_t^2)_t)) \\
&\qquad\qquad + \int_{\bar\CC(T_0)\times\bar\CC(T_0)} \big(D^{(0)}(\mu)\big)^2 \mu({\rm d}((c_t^1,w_t^1)_t,(c_t^2,w_t^2)_t))\bigg)^{1/2}.
\end{align*}
Since $\mu$ was arbitrary, we can take the infimum on both sides and arrive at
\begin{align*}
&\inf_{\mu\in\PP_{\zeta_{T_0}^\rho(\eta^1),\zeta_{T_0}^\rho(\eta^2)}(\bar\CC(T_0)\times\bar\CC(T_0))} \bigg(\int_{\bar\CC(T_0)\times\bar\CC(T_0)} \sup_{s\le T_0} \|(c_s^1,w_s^1)-(c_s^2,w_s^2)\|^2 \mu({\rm d}((c_t^1,w_t^1)_t,(c_t^2,w_t^2)_t))\bigg)^{1/2} \\
&\le e^{C T_0} \bigg(\big(T_0 C \DD_{T_0}(\eta^1,\eta^2) \big)^2 \\
&\qquad\qquad+ \inf_{\mu\in\PP_{\zeta_{T_0}^\rho(\eta^1),\zeta_{T_0}^\rho(\eta^2)}(\bar\CC(T_0)\times\bar\CC(T_0))} \Big(2 T_0 C \DD_{T_0}(\eta^1,\eta^2) \cdot\int_{\bar\CC(T_0)\times\bar\CC(T_0)} D^{(0)}(\mu) \mu({\rm d}((c_t^1,w_t^1)_t,(c_t^2,w_t^2)_t)) \\
&\qquad\qquad\qquad\qquad\qquad\qquad\qquad\qquad\qquad+ \int_{\bar\CC(T_0)\times\bar\CC(T_0)} \big(D^{(0)}(\mu)\big)^2 \mu({\rm d}((c_t^1,w_t^1)_t,(c_t^2,w_t^2)_t))\Big)\bigg)^{1/2}.
\end{align*}
Here we point out that $\zeta_{T_0}^\rho(\eta^1)$ and $\zeta_{T_0}^\rho(\eta^2)$ projected at time zero have the same distribution $\nu$.  We can rewrite both probability measures by disintegrating them with respect to the distribution of the initial value as it is done in the proof of \cite[Lemma 3.4]{budconroy} to show uniqueness. After combining both disintegrated conditional distributions under the same starting value, we obtain a coupling $\mu^*\in\PP_{\zeta_{T_0}^\rho(\eta^1), \zeta_{T_0}^\rho(\eta^2)}(\bar\CC(T_0)\times\bar\CC(T_0))$ such that for $((c_t,w_t)_t,(c'_t,w'_t)_t)\sim\mu^*$ it holds that $(c_0,w_0) = (c'_0,w'_0)$ $\mu^*$-almost surely. This proves that both integral terms vanish under the infimum and the assertion follows.
\end{proof}

\section*{Acknowledgments}
The authors thank the anonymous referees for their valuable comments and suggestions. In particular, it was noticed that the local martingale problem simplifies drastically in comparison to the Brownian setting, which allowed us to substantially shorten some arguments.
The authors thank P.\ Jung for very fruitful discussions in the early phase of the project. Further, DW acknowledges the financial support of the CogniGron research center and the Ubbo Emmius Funds (Univ.\ of Groningen). CH was supported by a research grant (VIL69126) from VILLUM FONDEN.

\end{document}